\documentclass[11pt]{article}

\usepackage{amsfonts, amsmath, amssymb,url, youngtab,bbm, color,blkarray, enumitem}
\usepackage[colorinlistoftodos]{todonotes}

\usepackage[left,pagewise,mathlines]{lineno}

\def\half{\tfrac{1}{2}}
\def\third{\tfrac{1}{3} }
\def\qed{ \hfill QED}

\newtheorem{thm}{Theorem}[section]
\newtheorem{prop}[thm]{Proposition}
\newtheorem{lemma}[thm]{Lemma}

\newtheorem{remark}[thm]{Remark}
\newtheorem{conj}[thm]{Conjecture}
\newtheorem{defn}[thm]{Definition}

\def\pmx{\begin{pmatrix}}
\def\emx{\end{pmatrix}}
\def \ssq{\begin{subequations}}
\def\esq{\end{subequations}}

\def\be{\begin{eqnarray}}
\def\ee{\end{eqnarray}}
\def\bee{\begin{eqnarray*}}
\def\eee{\end{eqnarray*}}
\def\pmx{\begin{pmatrix}}
\def\emx{\end{pmatrix}}

\def\bra{\langle}
\def\ket{\rangle}
\def\kb{ \ket \bra }
\def\nn{\nonumber}

\def\one{{\mathbbm 1}}
\def\eof{{\rm EoF}}

\def\id{{\mathcal I}}
 \def\1rt2{\tfrac{1}{\sqrt{2}}  }     

\newcommand{\proj}[1]{ | #1 \kb  #1|}
\newcommand{\norm}[1]{ \| #1  \|}
\def\tr{\hbox{Tr} \, }
\def\trp{\hbox{Tr} }
\def\ot{\otimes}
\def\op{\oplus}
\def\wh{\widehat}
\def\wtd{\widetilde}
\def\ovb{\overline}
\def\spn{{\rm span}}
\def\ds{\displaystyle}

\def\rmid{{\rm I }}
\def\bU{ {\bf U}}  
  
\def\bW{ {\bf W}}  

\title{ Extreme Points and Factorizability for  \\   New Classes of Unital Quantum Channels}

\author{Uffe Haagerup\footnote{Uffe Haagerup initiated and contributed significantly to this work until his untimely death  on July 5, 2015. He was Professor of Mathematics at the University of Copenhagen and the University of Southern Denmark.} \\  ~~   \\
\and Magdalena Musat  \footnote{ Department of Mathematical Sciences, University of Copenhagen, Universitetsparken 5, DK-2100, Copenhagen Ø, Denmark, musat@math.ku.dk} \\  ~~ \\
\and Mary Beth Ruskai \footnote{Department of Mathematics,   University of Vermont,   Burlington, Vermont  05405 USA,  \newline   mbruskai@gmail.com}}

\begin{document}

\maketitle

\begin{abstract}
We introduce and study two new classes of unital quantum channels.
The first class describes a 2-parameter family of channels given by completely positive  (CP)  maps
$M_3({\bf C})\mapsto M_3({\bf C})$ which are both unital and trace-preserving.
Almost every member of this family is factorizable and extreme in the
set of  CP maps which are both unital and  trace-preserving, but is not extreme in
either the set of unital CP maps or the set of trace-preserving CP maps.

We also study a large class of maps which generalize the Werner-Holevo channel
for $d = 3$  in the sense that they are defined in terms of partial isometries of 
rank $d-1$.  Moreover, we extend this to maps whose Kraus operators  have the
form $t \, \proj{e_j} \oplus V $ with $V  \in  M_{d-1} ({\bf C}) $  unitary and $t \in (-1,1)$.
We show that almost every map in this class is extreme in both the set of unital  CP maps
 and the set of trace-preserving CP maps.  We analyze in detail a
particularly interesting family which is extreme unless $t = \tfrac{-1}{d-1}$.  For $d = 3$,
this includes a pair of channels which have a dual factorization in the sense that they
can be obtained  by taking the partial trace over different subspaces
after using the same unitary conjugation in $M_3({\bf C}) \ot M_3({\bf C})$.

\end{abstract}

\bigskip % \pagebreak

\tableofcontents

  \pagebreak

\section{Introduction}

It is by now well-established that completely positive, trace-preserving (CPT) maps on
matrix algebras play an important role in quantum information theory because they
describe the effect of noise on a quantum system.   The set of CPT maps $M_d({\bf C})\mapsto M_d({\bf C})$
is convex, as is its dual, the set of unital completely positive (UCP) maps, and their intersection,
for which we use the acronym   UCPT maps.

By the Choi-Jamiolkowski isomorphism the set of UCPT maps from $M_d({\bf C}) $ to $M_d({\bf C})$ is in 
one-to-one correspondence with the set of bipartite states on ${\bf C}_d \ot {\bf C}_d$ for which
both quantum marginals (also called {\it reduced density matrices}) are given by the maximally mixed state, $\tfrac{1}{d} \rmid_d$.
From an operator algebra perspective, every density matrix, i.e., every positive  operator
$\rho \in M_d({\bf C})$ generates a state $\phi$ on the operator algebra $M_d({\bf C})$ given by $\phi(A) = \tr A \rho$.
The maximally mixed state $ \tfrac{1}{d} \rmid_d $ corresponds to the ``tracial state''  $\tau (A) = \tfrac{1}{d} \tr A $ (where  
$\tau$ denotes the trace  normalized  so that $\tau(\rmid_d) = 1 $). 

For qubits, i.e., two-dimensional systems,
  K\"{u}mmerer \cite{K}  showed that the extreme points of the set  of UCPT maps  on $M_2({\bf C})$ are
 precisely those that correspond to unitary conjugations. For $d \geq 3$,  the convex structure of the  
 UCPT maps  is much more complex.
However, relatively little was known about the extreme points other than unitary conjugations.

Despite some evidence \cite{MW} that  there are extreme points of the
UCPT maps which are not extreme for either the CPT or UCP maps, no explicit
examples were given in the literature until   Ohno \cite{Ohno} presented two examples
which he attributed to Arveson.
   We present a new family of UCPT maps on $M_3({\bf C})$, 
parameterized  by a pair of elements of the unit ball on ${\bf C}_2 $.  We show that with a few exceptions, 
each of these maps  
is extreme in the set of UCPT maps on $M_3({\bf C})$, but not extreme for either the UCP
or CPT maps.  Moreover, we show that all maps in this family are factorizable  
\cite{HM2011, HM2015}.  When mapped to bipartite states on ${\bf C}_3 \ot {\bf C}_3$
via the Choi matrix,
both quantum marginals are $\tfrac{1}{3} \rmid_3 $ for which the von Neumann entropy
is $\log 3$.  We give upper bounds on the entanglement of formation (EoF) demonstrating 
that, although these states are entangled,  they are not  maximally
entangled states on ${\bf C}_3 \ot {\bf C}_3$.   %  Indeed, the EoF is  less than $ \log 2 < \log 3$.
 
 In a complementary direction we study a class of extreme points of the set of UCPT maps 
 for $d \geq 3$ whose Choi-Jamiolkowski matrix has rank $d$ (sometimes known as Choi rank).  
 These maps are not given by unitary conjugations  and are, in general, not factorizable.
 Following a suggestion in \cite[Section 2]{Rusk}, we consider partial isometries on $M_d({\bf C})$
 constructed from unitary matrices in $M_{d-1}({\bf C})$ and extend this construction in a
 natural way to a much larger class of maps whose Kraus operators  have the
form $t \, \proj{e_j} \oplus V $ with $V  \in  M_{d-1} ({\bf C}) $  unitary and $t \in (-1,1)$.

 We show that almost all such maps are extreme points for both the  UCP and CPT maps. 
 We analyze in detail maps constructed from the unitary operator $ 2 \proj{\one_{d-1}} - \rmid_{d-1} $, 
 where $| \one_d \ket $ is the vector whose elements are all $d^{-1/2} $, and show that they
 are extreme in both the set of CPT and UCP maps unless $t = \tfrac{-1}{d-1}$.  
 When $d $ is odd, we construct a quite different
 family from rank two permutations in $M_{d-1}({\bf C})$;  these maps
 are extreme for all $t \in (-1, 1) $ when $d > 3$.

One of several equivalent ways (described in    \cite[Appendix A]{KMNR})  of formulating the Stinespring representation for a CPT map 
$\Phi : M_{d_A}({\bf C}) \mapsto M_{d_B}({\bf C})$
  on matrix algebras uses an auxiliary space ${\cal H}_E = {\bf C}_{d_E}$   called the environment and a
  unitary matrix  $U $ in $M_{d_B}({\bf C}) \ot M_{d_E}({\bf C})$ such that 
  \be   \label{stine}
       \Phi(\rho) =  ( \id \ot \trp_E )  U^* (  \rho \ot \proj{ \phi_E } ) U .
  \ee
  This is a natural  model for noise when  $d_A = d_B $,     the system
  is initially in a pure product state  $| \psi_A \ot \phi_E \ket $ on $ {\bf C}_{d_A}  \ot {\bf C}_{d_E} $,
  $U(t) $ describes the time evolution of the interacting system and environment, and $U = U(t_f) $
  corresponds to some later time $t_f$.
The question of factorizability of a UCPT map  (over a matrix algebra) roughly asks if the ancilla state $\proj{ \phi_E }$ can be
replaced by the maximally mixed state $\tfrac{1}{d_E} \rmid_{d_E} $.  When $ d_A = d_E$,
this is a natural model for  noise when the system  is initially decoupled from the environment which is in a 
maximally mixed or ``thermal'' state.  See, e.g., \cite{HHO, MKZ,SM}.   

The concept of a factorizable map was introduced by Anantharaman-Delaroche in \cite{A-D}  in a more general mathematical setting. 
It was further studied extensively in the 
context of unital quantum channels on $M_d(\mathbf{C})$ in  \cite{HM2011} and \cite{HM2015}. 
Following  \cite{HM2011}, we say that a UCPT map $\Phi$ on $M_d(\mathbf{C})$ has an exact factorization through $M_d(\mathbf{C}) \otimes {\mathcal{N}}$, 
where $\mathcal{N}$ is a von Neumann algebra with a normalized faithful trace $\tau$ (so that $\tau(I_{\mathcal{N}}) = 1$), 
if there is a unitary $\mathbf{U} \in M_d(\mathbf{C}) \otimes \mathcal{N}$ such that for all $\rho \in M_d(\mathbf{C})$,
\be    \label{exactfact}
\Phi(\rho) = (\mathcal{I} \otimes \tau)\mathbf{U}^*(\rho \otimes I_{\mathcal{N}}) {\mathbf{U}}.
\ee  
When $\mathcal{N} = M_\nu(\mathbf{C})$ is a matrix algebra, 
this is equivalent to $\Phi(\rho) = (\mathcal{I} \otimes \mathrm{Tr})\mathbf{U}^*(\rho \otimes \frac{1}{\nu}  \rmid_\nu) {\mathbf{U}}$.
 In \cite{MR} it was shown that for every $d \ge 11$, there are UCPT maps with an exact factorization through $M_d(\mathbf{C}) \otimes \mathcal{N}$, 
 where $\mathcal{N}$ is a von Neumann algebra of type II$_1$ which cannot be replaced by any  finite dimensional von Neumann algebra. 
 In view of \cite[Theorem 3.7]{HM2011} and the  announcement in \cite{JMVWY} that the Connes Embedding Problem has a negative answer,
  there are factorizable maps for some (presumably very large) $d$ which cannot even be approximated by maps with exact factorizations through
  matrix algebras.

  Most of the maps considered in Section~\ref{sect:part-isom}
 are extreme in the set of UCP or CPT maps and, hence, not factorizable.  However,  our work led us to
 a pair of channels,
 which have exact factorizations  through $M_3({\bf C} )\ot  M_3({\bf C} )$ that are dual 
 in the sense that they can be obtained   from the same unitary operator
 by  exchanging the roles of the system and auxiliary spaces.

  This paper is organized as follows.  In Section~\ref{sect:extUCPT}, we describe and study a
  family of UCPT maps which are not extreme in either the set of UCP maps or the set of CPT maps,
  but are extreme in the set of UCPT maps.   In Section~\ref{sect:over} we introduce
  UCPT maps $\Phi(\rho)  = \tfrac{1}{d - 1 + t^2 }  \sum_{m = 1}^d A_m^* \rho A_m $ with 
$A_m = t \, \proj{e_j} \oplus V_m $,  with $V_m  \in  M_{d-1} ({\bf C}) $  unitary, and $t \in (-1,1)$.
In Section~\ref{sect:dense} we prove several theorems which imply that, in general, 
maps of this form are extreme in both the set of UCP and the set of CPT maps
using the equivalent condition from \cite[Theorem 5]{Choi} of linear independence of the sets  $\{ A_m^* A_n \} $ 
and  $\{ A_m A_n^* \} $, respectively.
In Section~\ref{sect:crit} we consider explicit examples and subclasses of maps of this type. 
Section~\ref{sect:key} is devoted to the special case in which 
$V_m = 2 \proj{\one_{d-1} } - \rmid_{d-1} $.  Sections~\ref{sect:structure} to \ref{sect:espace} present the
details of our analysis for this case when $d > 3$ and $t \neq \tfrac{-1}{d-1} $; Section~\ref{sect:d=3} deals with the 
  case $d = 3$; and Section~\ref{sect:spec} deals with $t  = \tfrac{-1}{d-1} $ when $d > 3$.
  
  Appendix~\ref{app:factcond} discusses
  exact factorizations; Appendix~\ref{app:dual}
  considers dual pairs of channels $\Phi, \Psi $ associated with a unitary $\bU \in  M_p({\bf C}) \ot M_q({\bf C}) $;
  Appendix~\ref{app:choi4} extends a result in \cite{HM2011} to show that if a channel has Choi rank $\leq 4$, then both
  $\Phi \circ \Phi^* $ and $\Phi^* \circ \Phi $ are factorizable.   Finally, Appendix~\ref{app:B} 
  presents an example to show that  a set  $\{ A_m^* A_n \} $  can be linearly independent
when $\{ A_m A_n^* \} $  is linearly dependent.
 %  \pagebreak

\section{High rank extreme points of unital quantum channels} \label{sect:extUCPT}

\subsection{Background}

It is a fundamental result of Choi \cite{Choi}  that  when a CP map  $\Phi: M_{d_1}({\bf C}) \mapsto M_{d_2} ({\bf C})$ 
is written in the form  $\Phi(\rho) = \sum_k A_k^* \rho A_k $, then it is extreme in the  
convex set of CP maps for which $\Phi(\rmid_{d_1} ) = B $ for some fixed $B \in M_{d_2} $ if
the operators $\{ A_k \} $   satisfy  $\sum_k A_k^* A_k = B $ and  can be chosen so that the set
$\{ A_m^* A_n \} $ is linearly independent.   Thus, a UCP map is extreme if and only if $\{ A_m^* A_n \} $
 is linearly independent, and a CPT map is extreme if and only if $\{ A_n A_m ^*\} $ is linearly
independent.   This implies that a  UCPT map whose Choi-rank is $> d_1$ can not be an extreme point
of the UCP maps and one whose Choi-rank is $> d_2$ can not be an extreme point of the CPT maps.
Thus a UCPT map $\Phi: M_d({\bf C})\mapsto M_d({\bf C})$  with Choi rank $> d $ can not be an extreme point of
either the UCP or CPT maps.

Nevertheless  this does not preclude the possibility that a UCPT map with Choi-rank   greater than $d$ can be
  extreme in the set of UCPT.  We present a family of such maps for $d = 3$ 
 with Choi-rank   equal to $4$, which are
  parameterized by  $\alpha, \beta \in {\bf C} $ with $ |\alpha|^2 +| \beta|^2  = 1$.
  When  $|\alpha|^2  \neq  0, \half, 1$ these maps are
  extreme in the set of UCPT maps.  Moreover, all of the maps in this family have
   exact factorizations through  $M_3({\bf C})  \ot M_2({\bf C}) $.
  
  Although earlier work \cite{MW} suggested the existence of such maps, the only explicit examples in the
  literature are due to Arveson, as presented by Ohno in \cite{Ohno} for $d = 3 $ and $d = 4$.  For comparison
  with our results, 
  we note that in the $d = 3$ case the Choi-Kraus operators for  the Arveson-Ohno map are
   \begin{align}   \label{ohnomap}
     A_1 & = \proj{e_1} & A_2 & = |e_1 \kb e_2 | + \sqrt{2} \, |e_2 \kb e_3|   \nn \\
     A_3 & = \sqrt{2} \, |e_2 \kb e_1 | + \sqrt{3} |e_3 \kb e_2 | & A_4 & = |e_3 \kb e_1 |+ \sqrt{2} \, |e_1 \kb e_3 |
 \end{align}
In  Section~\ref{sect:remarks} we show that the map $\Phi(\rho) = \sum_{k=1}^4   A_k^* \rho A_k $ is {\em not} factorizable, i.e.,
it does not have an exact factorization through  $M_3({\bf C})  \ot  {\cal N} $ for any von Neumann algebra ${\cal N} $.
  However, it follows from Proposition~\ref{prop:choi4}, which is a straightforward generalization of  \cite[Remark~5.6]{HM2011}, 
   that the maps $\Phi \circ \Phi^*$and $\Phi^* \circ \Phi$   
   have  exact factorizations through  $M_3({\bf C})  \ot M_4({\bf C}) $.
    
 \subsection{A factorizable family of extreme UCPT maps}   \label{sect:UCPT}
 
 We present a new family of UCPT maps $M_3({\bf C}) \mapsto M_3({\bf C})$ which have Choi-rank equal to $4$ so that
 they can not be    extreme in either the set of UCP or CPT maps.  However, these maps are extreme in the set of UCPT maps.
All  maps in this family, including those which are not extreme, have  exact factorizations through  $M_3({\bf C})  \ot M_2({\bf C}) $.
 
 \begin{thm}  \label{thm:UCPT}
 Let    $\alpha, \beta \in {\bf C}$   with $ | \alpha |^2  +  | \beta|^2 = 1 $ and let
 \begin{align}   \label{UCPT}
    A_1 & = \alpha \proj{e_1} + ~ |e_2 \kb e_3 |  &    A_2 & = \beta |e_1 \kb e_3 | + ~ |e_3 \kb e_2 |   \nn  \\
   A_3& =  -   |e_1 \kb e_2 |  -  \ovb{ \beta} | e_3 \kb e_1 |  & A_4 &  = ~  |e_2 \kb e_1 | + \ovb{\alpha } \proj{ e_3 } ~.
\end{align} 
Then   $\Phi_{\alpha,\beta}(\rho) =  \half \ds{ \sum_{k = 1}^4 A_k^* \rho A_k }$  
  is an extreme point in the set of UCPT maps if $|\alpha|^2  \neq  0, \half, 1$.
 Moreover, $\Phi_{\alpha,\beta}$ has an exact factorization through  $M_3({\bf C})  \ot M_2({\bf C}) $
for all $\alpha, \beta $.
 \end{thm}   
  \noindent{\bf Proof:}  To show that  $\Phi_{\alpha,\beta} $ has an exact factorization through  $M_3({\bf C})  \ot M_2({\bf C}) $ let
 \be   \label{unitaryUCPT}
 {\bf U} =\sum_{j,k \in 1,2}  A_{2(j-1)+k} \ot |e_j \kb e_k | =   \pmx  A_1& A_2  \\  A_3 & A_4  \emx  .
 \ee
Then  $  {\bf U}  \in M_3({\bf C}) \ot M_2({\bf C}) \simeq M_6({\bf C}) $ is unitary and 
  \be   \label{fact:UCPT}
      \Phi_{\alpha,\beta} (\rho) = \half \sum_{k = 1}^4 A_m^* \rho A_m = (\id_3 \ot \tr )( U^* \big(\rho \ot \half \rmid_2 ) U \big).
  \ee

  To  prove the rest  of this theorem, recall that Landau and Streater \cite{LS} showed that a 
  necessary and sufficient conditions for $\Phi_{\alpha,\beta} $ to be
  an extreme UCPT map
 is that  the set $\{  B_{jk}  = A_j^* A_k \op A_k A_j^* \} $ is linearly independent in $M_3({\bf C}) \op M_3({\bf C}) $.
  First, observe that the diagonal of $B_{jk} $  is zero  if $j \neq k$, and that each $B_{jj}$ is diagonal
  with diagonal elements corresponding to the rows of the matrix
  \bee  \pmx
     |\alpha|^2 & 0 & 1 &   |\alpha|^2  & 1 & 0 \\ 0 & 1 & |\beta|^2 &  |\beta|^2 & 0 & 1 \\
       |\beta|^2  & 1 & 0 & 1 & 0 &  |\beta|^2  \\ 1 & 0 & |\alpha|^2 & 0  & 1 & |\alpha|^2
\emx   ~.\eee
     After deleting the third and fourth columns one obtains a matrix whose determinant is $- |\alpha|^2 |\beta|^2 $.
    This implies that if $|\alpha | \neq 0, 1 $ the rows are linearly independent and, hence, the set $\{ B_{jj} \}_{j=1}^4 $ is linearly independent.
    
  Next observe that 
  \be 
      B_{12}   =   A_1^* A_2 \op A_2 A_1^*  & = &     \ovb{\alpha} \beta  \, |e_1 \kb e_3 | \op 0_3    + \beta\, 0_3 \op |e_1 \kb e_2 |    \nn  \\
    B_{41}   =  A_4^* A_1 \op A_1 A_4^*& = &      \quad ~ \,  |e_1 \kb e_3 | \op 0_3    + \alpha \,  0_3 \op |e_1 \kb e_2 |    +     \alpha \, 0_3 \op |e_2 \kb e_3  |\\  \nn
  B_{34}   =   A_3^* A_4 \op A_4 A_3^* & = &  -  \ovb{ \alpha } \ovb{\beta} \, |e_1 \kb e_3 | \op 0_3  \qquad  \qquad  \qquad  \quad  ~ - \ovb{\beta} \, 0_3  \op |e_2 \kb e_3 |
\ee
which are clearly  linearly independent if and only if 
$$\det \pmx   \ovb{\alpha} \beta &   \beta  & 0 \\ 1 & \alpha &  \alpha  \\   \ovb{\alpha}  \ovb{\beta} & 0 & \ovb{\beta} \emx =   |\beta|^2 (  2 |\alpha|^2  - 1) \neq 0 $$
Since  $B_{jk} = B_{kj}^*$, the matrices $\{ B_{21} , B_{14}, B_{43} \} $ are also linearly independent under the same conditions.
Next, we similarly treat
\bee
  B_{24} =    A_2^* A_4 \op A_4 A_2^* & = &    \qquad \qquad  \qquad \qquad ~   \ovb{\alpha }\,  |e_2 \kb e_3 | \op 0_3    +      \ovb{\alpha }  \ovb{ \beta} \,  0_3 \op  |e_3 \kb e_1 |  \\
   B_{32} =  A_3^* A_2 \op A_2 A_3^* & = & -  \beta \, |e_1 \kb e_2 | \op 0_3  + \beta  \, | e_2 \kb e_3 | \op 0_3  ~ -  \quad ~ \, 0_3 \op |e_3 \kb e_1 | \\
   B_{13} =   A_1^* A_3 \op A_3 A_1^* & = &-  \ovb{\alpha}  \, |e_1 \kb e_2 | \op 0_3     \qquad  \qquad  \qquad  \quad  ~  -  \ovb{\alpha }  \ovb{ \beta}  \, 0_3 \op  |e_3 \kb e_1 | 
\eee
which are linearly independent if and only if
$$ \det \pmx  0 & \ovb{\alpha}   & \ovb{\alpha}    \ovb{\beta} \\ \beta & \beta & 1 \\   \ovb{\alpha} & 0 &   \ovb{\alpha} \ovb{\beta} \emx =  
   \ovb{\alpha}^2 ( 1 - 2 |\beta|^2)  \neq 0 ~.$$
After again observing that the adjoints are linearly independent  under the same conditions, we can conclude that  if $|\alpha|^2 \neq  0, 1, \half $  
then each of the four sets  
$$\{ B_{12}, B_{41}, B_{34} \}, \qquad  \{ B_{21}, B_{14}, B_{43} \} , \qquad   \{ B_{24}, B_{32}, B_{13} \} , \qquad \{ B_{42}, B_{23} , B_{31} \} $$
consists of linearly independent matrices in $M_6({\bf C})$.  Moreover, the only common point in the spans of each of these
four sets is the zero matrix $0_6 $.  Therefore $\{ A_j^* A_k \op A_k A_j^* \}_{j \neq k } $ is linearly independent when 
$|\alpha|^2 \neq 0, \half , 1 $.  Combining this with our observations above for $B_{kk}$, implies that
$\big\{ B_{jk} :  j , k \in \{ 1, 2, 3, 4  \} \big\} $ is linearly independent.   \hfill QED

Because the maps $\Phi_{\alpha,\beta} $ are parameterized by a unit vector in ${\bf C}_2 $, it is tempting to
associate each channel with a qubit state.  However, the vectors $ |v \ket $ and $ e^{i \theta} | v \ket $
represent the same physical state.  But the channels associated with, e.g., $\alpha = 1$ and $\alpha = -1 $
are not the same.

% \pagebreak

\subsection{Remarks} \label{sect:remarks}
\subsubsection{Non-extreme cases}  

When $|\alpha| = 1 $ or $|\beta|  = 1 $, the map $\Phi_{\alpha,\beta} $ can be written as a convex combination of unitaries, e.g.,
when $\alpha = 1$, $\Phi_{1,0} = \tfrac{1}{4} \sum_{k=1}^4 U_k^* \rho U_k $ where
\begin{align*}
%\beta & = 1   \quad  &  U_1 & = \pmx 0 & 0 & 1 \\ 1 & 0 & 0 \\ 0 & 1 & 0 \emx    & U_2 & = \pmx 0 & 0 & 1 \\ -1 & 0 & 0 \\ 0 & 1 & 0 \emx
%  & U_3 & = \pmx 0 & 1 & 0 \\ 0 & 0 & 1 \\ 1 & 0 & 0 \emx & U_4 & = \pmx 0 & 1 & 0 \\ 0 & 0 & -1 \\ 1 & 0 & 0 \emx,  \\  ~~ \\
    U_1 & = \pmx 0 & 1 & 0 \\ 1 & 0 & 0 \\ 0 & 0 & 1\emx 
     &  U_2 & = \pmx 0 & -1  & 0 \\    1 & 0 & 0 \\ 0 & 0 & 1\emx      &  U_3 & = \pmx 1 & 0 & 0 \\ 0 & 0 & 1 \\ 0 & 1 & 0 \emx
        &  U_4 & = \pmx  1 & 0 & 0 \\ 0 & 0 & 1  \\  0 & -1 & 0 \emx.
\end{align*} 

However, when $\alpha = \beta = \1rt2$,   the map $\Phi_{\alpha,\beta} $ is not a convex combination of unitaries.  Let   $| \one_4 \ket $ denote the vector with
all elements $ \half $.  Then  $W = 2 \proj{\one_4} - \rmid_4 $ is unitary.  Now let
$X_j =\1rt2  \sum_{k =1 }^4 w_{jk} A_k $, then
\begin{align*}
   &  \quad     &     X_1 & =   \tfrac{1}{4} \pmx -1 & -\sqrt{2} & 1 \\ \sqrt{2} & 0 &-  \sqrt{2} \\ -1 & \sqrt{2} & 1\emx     &
     X_2 & =  \tfrac{1}{4}  \pmx 1 & -\sqrt{2}  &  -1  \\ \sqrt{2} & 0 & \sqrt{2}  \\ -1 & - \sqrt{2} & 1 \emx    \\
   ~~   \\  %   \alpha  = \beta  & = \1rt2  \\
 &  & X_3 & = \tfrac{1}{4} \pmx 1 &  \sqrt{2} & 1 \\ \sqrt{2} & 0 & \sqrt{2}  \\  1 &  \sqrt{2} & 1  \emx &
  X_4 & =   \tfrac{1}{4} \pmx 1 &  - \sqrt{2} & 1 \\  - \sqrt{2} & 0 & \sqrt{2}  \\ - 1 &  \sqrt{2} &- 1  \emx ~. 
\end{align*}
Since $ W$ is unitary, $  \sum_j X_j^* \rho X_j  = \half  \sum_k A_k^* \rho A_k = \Phi(\rho) $.  
Although $\wtd{X}_2 = 2  X_2 $ is unitary,  this is not true for $j = 1, 3, 4$.  In fact, both   $\{ X_j^* X_k \}_{j,k = 1,3,4} $ 
and  $\{ X_j X_k^* \}_{j,k = 1,3,4} $  are linearly independent sets
so that  $\Psi(\rho) = \tfrac{4}{3} \sum_{k = 1,3,4} X_j^* \rho X_j $ is an extreme point of both the UCP  and the CPT maps. 
 Thus, 
$\Phi(\rho) = \tfrac{1}{4} \wtd{X}_2^*(\rho) \wtd{X}_2   + \tfrac{3}{4} \Psi(\rho)$ is a convex combination of a unitary conjugation and a
UCPT map with Choi rank $3$.  (It may be worth noting that $ \Psi $ is an extreme point of both the UCP and CPT maps which has
Choi-rank $3$, but is not of the form considered in Section~\ref{sect:part-isom}.)

The argument above easily extends to  $\alpha =  \pm \ovb{\beta} $  with $2 X_3 $  unitary.
 It seems reasonable to conjecture that,
whenever  $|\alpha | = |\beta| = \1rt2$, the channel $\Phi_{\alpha,\beta} $ is a convex combination of a
unitary conjugation and a map with Choi-rank $3$ which is extreme in both the UCP and CPT maps.
However, one would need a different unitary transformation relating $\{ A_k \} $ to $\{ X_k \} $.

\subsubsection{Entanglement of Formation}

Recall that the Choi matrix associates any CP map with a bipartite state  $\rho$ on ${\bf C}_d \ot {\bf C}_d $ and
those  associated with UCPT maps have quantum marginals given by the
maximally mixed state $\frac{1}{d} \rmid_d $.   Since the UCPT maps considered here are not associated
with pure maximally entangled states, any measure of entanglement based on entropy
will be less than the maximum value of $\log d$.  A natural one to use in this situation is
the {\em entanglement of formation} \cite{eof}, defined as
%\todo[inline, color=yellow]{Normally one uses
%$\log_2 $ but it is really irrelevant for the purpose of making comparisons as long as one is consistent.}  
\be
    \eof(\rho) = \inf \bigg\{ \sum_k p_k E(\psi_k) :   \rho = \sum_k p_k \proj{\psi_k} \bigg\}
\ee
where  $E(\psi_k) = S\Big( \trp_B \proj{\psi_k} \Big) $ and $S(\rho) = - \tr \rho \log \rho $ is the von Neumann
entropy.  An upper bound on the EoF (which is probably optimal) for the Choi matrix of the UCPT maps 
 is readily calculated
for the examples given here since  the Kraus operators, which correspond to eigenvectors of the Choi matrix,  
are all associated with pure bipartite states. 
%  \todo[inline, color=pink]{I am not familiar with this. } 
%\todo[inline, color=yellow]{This is how Choi obtained the operator sum decomposition, i.e., by
% ``stacking'' the eigenvectors of $J_\Phi $ to obtain $A_k$ which can also be regarded as the 
% coefficients of a bipartite state in ``Schmidt form'', i.e.,  $\psi_k = \sum_{ij} (A_k)_{ij} |e_i  \ot f_j \ket $ 
% for any pair of O.N. bases.   When $J_\Phi $ is interpreted as a density matrix on ${\bf C}_d \ot {\bf C}_d $
%  these become the eigenvectors and one obtains the spectral decomposition after adjusting the normalization.}

For our family of channels, the Choi matrix is
\bee
  \rho = \frac{1 + |\alpha|^2}{6} \proj{\psi_1 } +  \frac{1 + |\beta|^2}{6} \proj{\psi_2 }  + \frac{1 + |\beta|^2}{6} \proj{\psi_3}  + \frac{1 + |\alpha|^2}{6} \proj{\psi_4 }, 
  \eee
  where 
  \begin{align*}
      \psi_1 & =  \frac{1}{\sqrt{1 + |\alpha|^2}} \big( \alpha \, |e_1 \ot e_1  \ket + |e_2 \ot e_3 \ket \big) 
        &  \psi_2 & =   \frac{1}{\sqrt{1 + |\beta|^2}} \big( \beta \, |e_1 \ot e_3  \ket + |e_3 \ot e_2 \ket \big) \\
            \psi_3 & =   \frac{1}{\sqrt{1 + |\beta|^2}} \big(  \, |e_1 \ot e_2  \ket +  \ovb{\beta} |e_3 \ot e_1\ket \big)  
           &     \psi_4 & =   \frac{1}{\sqrt{1 + |\alpha|^2}} \big(   |e_2 \ot e_1  \ket + \ovb{\alpha}  |e_3 \ot e_3 \ket \big)~.
  \end{align*} 
This implies that
\bee
     \eof(\rho_{AB})  & \leq &  
            \tfrac{ 1 + |\alpha|^2 }{3}   h \big(\tfrac{1}{1 + |\alpha|^2} \big)  +  \tfrac{ 1 + |\beta|^2 }{3}   h \big(\tfrac{1}{1 + |\beta|^2} \big) 
\eee
where $h(x) = - x \log x - (1-x) \log (1-x) $ is the binary entropy.
For $|\alpha|^2 = |\beta|^2 = \half $ the expression above takes its largest value of $h\big(\frac{1}{3} \big) = 0.918296 <  1 = \log 2$, and
for $|\alpha|^2  = 0, 1 $ its smallest value of $\tfrac{2}{3} h\big( \half \big) = \tfrac{2}{3} \approx 0.66667.$
This is considerably less than the EoF of a maximally entangled state which is $\log 3 \approx 1.58496$.  
 Thus, we have a large  family of bipartite states on ${\bf C}_3 \ot {\bf C}_3 $ 
which are extreme in the convex set of states whose quantum marginals 
 are $\tfrac{1}{3} \rmid_3 $, but whose entanglement as measured by the EoF 
can be less than half that of a maximally entangled state.

For comparison, we can bound the EoF for  the Arveson-Ohno example \eqref{ohnomap} by similarly observing that
its  Choi matrix is
\bee
     \rho = \tfrac{1}{12} \proj{e_1 \ot e_1 } +   \tfrac{1}{4} \proj{\psi_2} +  \tfrac{5}{12} \proj{\psi_3} +  \tfrac{1}{4} \proj{\psi_4} ,
\eee
with $ \quad     \psi_2 = \tfrac{1}{\sqrt{3} }\big( |e_1 \ot e_2 \ket + \sqrt{2} \, |e_2 \ot e_3 \ket \big), \quad
   \psi_3   =    \tfrac{1}{\sqrt{5} }\big(\sqrt{2} \,  |e_1 \ot e_2 \ket + \sqrt{3} \, |e_3 \ot e_2 \ket \big)$,  and \linebreak
      $    \psi_4=  \tfrac{1}{\sqrt{3}} \big( |e_3 \ot e_1\ket + \sqrt{2} \, |e_1 \ot e_3 \ket \big)$. 
This implies  that
\bee
   \eof(\rho_{AB}) \leq \tfrac{1}{4} h \big(\tfrac{1}{3} \big) + \tfrac{5}{12} h \big(\tfrac{2}{5} \big) + \tfrac{1}{4} h \big(\tfrac{1}{3} \big)  \approx 0.8637.
\eee

\subsubsection{The Arveson-Ohno channel is not factorizable}

\begin{prop}
The UCPT channel $\Phi(\rho) = \sum_{k=1}^4  A_k^* \rho A_k $   with $A_k$ given by \eqref{ohnomap}
is not factorizable.
\end{prop} 
 \noindent{\bf Proof:} 
  First observe that the matrices $\{ A_1, A_2, A_3, A_4 \} $ are linearly independent,
so that we can use Proposition~\ref{prop:factcond}  in Appendix~\ref{app:factcond}.
Let  $ {\cal N} $ be a von Neumann algebra with normal faithful trace satisfying $\tau(\rmid_{\cal N} ) = 1$.
We  assume that we can find a set of matrices  $\{ Y_1, Y_2, Y_3, Y_4 \} $ in  ${ \cal N}$ which satisfy $\tau( Y_j^* Y_k ) =   \delta_{jk} $, 
and  find a contradiction.
 A straightforward calculation shows that  $\bra e_3,  A_j^* A_k \, e_2  \ket  = 0$ unless $j = 4, k = 2 $
          which    implies  $ 0 = \1rt2  Y_4^* Y_2 $.  Multiplying this by  $Y_4$ on the left and $Y_2^* $ on the
           right implies $Y_4 Y_4^* Y_2 Y_ 2^* = 0 $.
           
           Next, observe that $ \bra e_m,  A_j^* A_k e_m \ket = 0 $ if $j  \neq  k $, and
  \bee
  \bra e_2 , A_j A_j^* \, e_2  \ket = \begin{cases}    \half &  j =2  \\ \half & j = 3  \\  0 &  j = 1, 4 \end{cases} \qquad  \qquad
   \bra e_3 , A_j A_j^* \,  e_ 3  \ket = \begin{cases}   \tfrac{3}{4} &  j = 3  \\  \tfrac{1}{4}  &  j = 4  \\  0 & j = 1,2\end{cases} 
   \eee              
  from which it follows   that
  \be
       \half Y_2 Y_2^* + \half Y_3 Y_3 ^* =  \rmid_{\cal N}   \qquad     \tfrac{3}{4} Y_3 Y_3^* + \tfrac{1}{4} Y_4 Y_4^* = \rmid_{\cal N}.
  \ee
Combining these to eliminate  $Y_3 Y_3 ^* $  gives $ 3Y_2 Y_2^* - Y_4 Y_4^* = 2\rmid_{\cal N} $.  Then multiplying by $ Y_4 Y_4^* $
and using   $Y_4 Y_4^* Y_2 Y_ 2^* = 0 $ implies    $ Y_4 Y_4^* = - (Y_4 Y_4^* )^2 $ which implies $Y_4 = 0$.  Thus, we
have shown  $Y_2 Y_2^* = \tfrac{2}{3} \rmid_{\cal N}  $.

Finally, we observe that $ \bra e_3,  A_j^* A_k e_3 \ket = 0 $ if $j  \neq  k $ and  
$  \bra e_3 , A_j^* A_j \,  e_ 3  \ket = \begin{cases}  \half &  j  = 2, 4 \\ 0 & j = 1,3 \end{cases}  $ which implies
$ \half Y_2^*  Y_2 +  \half  Y_4^* Y_4 = \rmid_{\cal N} $.  However, since $Y_4 = 0 $, this implies
$  Y_2^*  Y_2 = 2 \rmid_{\cal N} $ which is not consistent with $Y_2 Y_2^* = \tfrac{2}{3}\rmid_{\cal N}  $.  
 (Since, e.g., it would give  $\tau( Y_2 Y_2^* ) \neq \tau( Y_2^* Y_2 )$.)  \qed

\section{Extreme points from partial isometries}  \label{sect:part-isom}

\subsection{Definitions and basic properties }   \label{sect:over}

It is well-known that the extreme points of the convex set of UCPT maps include conjugation with a
single unitary, i.e., $\rho \mapsto U^* \rho U$, and that these are extreme in both the set of UCP maps and the set of CPT maps.
It is also known that there are maps which are not unitary conjugation but, nonetheless, extreme in both the
set of UCP and the set of CPT maps.   The simplest  example is the Werner-Holevo channel  \cite{WH} for $d = 3$ and its
symmetric counterpart.

One interesting class of quantum channels, which can be regarded as a generalization of this, was
proposed in \cite[Section 4.2]{Rusk}, but not widely studied.   These channels have
  $d$ Kraus operators, each of which is a multiple of a partial isometry of rank $d - 1$.  To define them, 
we will use the cyclic shift operator $S = \sum_k | e_k \kb e_{k+1} |$.  Let $\{ V_1, V_2, \ldots V_d \} $ be 
a set of unitary matrices in $M_{d-1} ({\bf C}) $, and let  
\be  \label{zeroform}
  A_m =    S^{-m+1} \pmx  0  & 0 \\ 0 & V_m \emx S^{m-1} , \qquad  m = 1, 2, \ldots d.
  \ee
Then it is easy to check that $A_m^* A_m = A_m A_m^*  =   ( \rmid_d -  \proj{e_m}) $, so that
$\sum_m  A_m^* A_m = \sum_m A_m A_m^*  = (d-1)  \rmid_d $ and the map
 $ \Phi(\rho) =   \tfrac{1}{d-1} \sum_m A_m^* \rho A_m $ is  both CPT and UCP.
  One can, more generally, consider operators of the form
\be  \label{tform}
  A_m =  S^{-m+1} \pmx  t  & 0 \\ 0 & V_m \emx S^{m-1} , \qquad  m = 1, 2, \ldots d
       \ee
with $ t \in [-1,1] $.  Since  $A_m^* A_m = A_m A_m^*  =   \rmid_d -  (1-t^2)  \proj{e_m} $, 
 $$\sum_m  A_m^* A_m = \sum_m A_m A_m^*  = (d - 1 + t^2 ) \rmid_d $$
 which implies that the map $\Phi$ given by
  \be   \label{tchan} 
      \Phi(\rho) =   \tfrac{1}{d - 1 + t^2 }  \sum_{m = 1}^d A_m^* \rho A_m 
       \ee
   is both UCP and CPT.  Note that the Kraus operators for the channel are
  $\wtd{A}_m \equiv \tfrac{1}{\sqrt{d-1 +t^2} } A_m $  which satisfy 
   $\sum_m \wtd{A}_m^*  \wtd{A}_m = \sum_m \wtd{A}_m  \wtd{A}_m^*  = \rmid_d $.

 \begin{remark}  \label{rmk:factpm1}
 When $t = \pm 1$, the matrices $ A_m $ in \eqref{tform} are unitary and the channel  \eqref{tchan}  can not be 
 extreme in either the UCP or CPT maps.    On the contrary, it is factorizable.  When
 ${\bf U} = \oplus_{m = 1}^d  A_m = \sum_{m=1}^d  A_m \ot \proj{e_m} $ is a block diagonal   $d^2 \times d^2 $ unitary matrix,
 \be   \label{factpm1}
      \Phi(\rho) =   \sum_m A_m^* \rho A_m  =  (\id_d \ot {\rm Tr} ) \big( {\bf U}^* (\rho \ot \tfrac{1}{d} \rmid  ){\bf U} \big).
 \ee 
 \end{remark}

We now focus on the case  $t \in (-1,1)$, for which we  are interested in the  linear independence of the sets
    $ \{ A_m^* A_n \} $  and    $ \{ A_m A_n^* \} $,  as these are precisely the conditions for a map
   $ \Phi $ given by \eqref{tchan} to be an extreme point of the set of UCP and CPT maps respectively.   
   
    {\begin{prop} \label{prop:diag} Let  $t \in (-1,1)$ and for  $m = 1, 2, \ldots d $ let  $\{A_m \}$  be as in (9) 
 with $\{V_m \}$  arbitrary  unitary matrices in ${ M}_{d-1}({\bf C})$. Then each $A_m$ 
  is a normal matrix in $M_d({\bf C})$, and $\{ A_m^*A_m \}_{m=1}^d $ is a  
  linearly independent set of diagonal matrices.  Moreover it is a basis for the set of diagonal matrices in $M_d({\bf C})$.
\end{prop}
\noindent {\bf Proof:} First, observe that
\bee  
A_m^*A_m =  A_m A_m^*  =     S^{-m+1} \pmx  t^2  & 0 \\ 0 & \rmid_{d-1}  \emx   S^{m-1}   
=\rmid_d - (1-t^2) \, \proj{e_m} ,
 \eee
 which implies that $\sum_{m=1}^d A_m^* A_m =  (d-1+t^2) \rmid_d$.   Thus, for the purpose of 
 determining  linear independence, one can replace   $A_m^* A_m $ by
$     A_m A_m^* - \rmid_d = - (1-t^2) \proj{e_m }$.
This immediately implies the desired linear independence.  \qed

\begin{remark} \label{rmk:diag}  It follows immediately  from the proof of Proposition~\ref{prop:diag}  that when $t \in (-1,1)$
\be
   \spn \{ A_m A_m^* \}_{m=1}^d  = \spn \{ A_m^* A_m  \}_{m=1}^d = \spn \{ \proj{e_m } \}_{m =1 }^d .
\ee    
Therefore, for the purpose of determining the linear independence of the set $\{ A_m^* A_n \} $,
one can make arbitrary changes to the diagonal elements of  $A_m^* A_n $ and it suffices to show that
the set $\{ A_m^* A_n \} _{m \neq n}$ is linearly independent.
 \end{remark}

 We study here only maps generated by unitaries in $M_{d-1}({\bf C})$.  If one allows
more general partial isometries, many more examples of maps which are 
extreme in both the set of UCP and CPT maps can be constructed.  For $d = 6$,
one could  construct three Kraus operators
\bee
      A_1 =  \half \pmx X & Y &  0_2 \\ Y & X & 0_2 \\ 0_2 & 0_2 & 0_2 \emx \qquad
         A_2 =\half  \pmx X &   0_2 & Y \\ 0_2 & 0_2 & 0_2 \\ Y &  0_2 & X \\  \emx 
         \qquad A_3  =  \half \pmx 0_2 & 0_2 & 0_2 \\ 0_2  & X & Y  \\0_2 &  Y & X   \emx 
\eee
based on the unitary  $ \half  \pmx X & Y \\ Y & X \emx$ with
$X = \pmx -1 & 1 \\  1 & -1 \emx, ~ Y = \pmx 1 & 1 \\ 1 & 1 \emx $ and 
related by  $A_m =  S^{-2} A_{m-1} S^2 = S^{-2(m-1)} A_1  S^{2(m-1)} $.

We note that Ohno \cite{Ohno} considered the channel with 
$$A_1 = \sqrt{\tfrac{d-2}{d-1} } \big(\rmid_d - \proj{e_1} \big), \qquad 
 A_k = \tfrac{1}{\sqrt{d-1}} \big( |e_1 \kb e_j | + | e_j \kb e_1 \big),  ~~ k = 2, \ldots d, $$
 which also has Choi rank $d$ and is constructed from partial isometries, albeit 
 with lower rank than $d-1$.

 \subsection{Extreme Points are Generic}   \label{sect:dense}
 
 \subsubsection{Overview}
 
 In this section we show that channels generated using \eqref{tchan} and \eqref{tform} are  usually extreme 
 in both the UCP and CPT maps in two ways.  
 First, we consider the  $V_m$ fixed and show that, in general,  for
 all but a finite number of $t \in (-1,1) $ the corresponding maps are extreme.  
 Then, we fix   $d$ and $t$ and show that almost all choices of $\{ V_m \} $ generate maps which are extreme.
 
 The latter implies that, roughly speaking, if we can find one choice of $V_k$ in  \eqref{tform}  for which
 the channel given by \eqref{tchan} is extreme in the UCP or CPT maps, then almost every choice of $V_k$ 
 also generates a channel which is extreme.  However (as discussed after Theorem~\ref{thm:density})
 if the example has all $V_j \neq V_k $, this does
 not preclude the  possibility that no channels with, e.g., $V_1 = V_2$ are extreme.   But an example
 with $V_1 = V_2 = \ldots = V_d $ which is extreme implies that almost every channel with unconstrained $V_j $
 is also extreme.  In Section~\ref{subsect:key}  we present a simple  example of this type which is extreme unless
 $t = \tfrac{-1}{d-1} $ and for which we also have $V_j = V_j^* $.
 In Section~\ref{subsect:rankone} we provide additional evidence for generic extremality even when  $t = \tfrac{-1}{d-1} $.

 \begin{remark}   \label{rmk:prelim}
 Before presenting these results precisely, we make some observations
 about maps on $M_d({\bf C}) $
 that will be used in both settings. 
  \end{remark}
\begin{enumerate}[label=\Alph{*})]
     \item   A map 
      $\Phi$ is extreme in the set of UCP maps if and only if its adjoint 
       $\Phi^* $  is extreme in the set of CPT maps.   Thus,  a UCPT map with  $\Phi = \Phi^* $ is extreme 
       in the UCP maps if and only if it is extreme in the CPT maps.
                
       \item    It follows from   \cite[Theorem 5]{Choi} that  a map
      $\Phi(\rho)  = \sum_k A_k^* \rho  A_k$ with $ \sum_k A_k^* A_k = \rmid_d $ 
       is extreme in the set of UCP maps if and only if the set of matrices $\{ A_m^* A_n \} $
        is linearly independent.       
       
        \item   It then follows from   \cite[Corollary~2.3]{ HM2011}  that if a UCPT map is extreme
        in either the UCP or  CPT maps, it is not factorizable.
       
       \item     A set of vectors $\{ v_j \} $ is linearly independent if and only if the Gram matrix 
       with elements $g_{jk} = \bra v_j , v_k \ket $ is non-singular.  Hence,
        the set of matrices  $\{ A_m^* A_n \} $ is linearly independent if and only if its
       Gram matrix  $G$  formed using the Hilbert-Schmidt inner product 
       with elements  $g_{jk,mn} = \tr A_j^* A_k (A_m^* A_n )^* $
       is non-singular.
  
 \end{enumerate}
 
 In view of (B) above, we present our results in terms of the linear independence of the sets
 $\{ A_m^* A_n \} $ and $\{ A_m A_n^*  \} $.  The implications for extremality of UCP and CPT
maps is straightforward.

 When  $V_m \neq V_m^*$ it is possible that $\{ A_m^* A_n \} $ is 
 a linearly independent set,  but $\{ A_m A_n^*  \} $ is linearly dependent. 
 Therefore, a UCPT map $\Phi$ can be extreme in the set of CPT maps,
 but not in the set of UCP maps (or vice versa).    An explicit example for $d = 4$ is
 given in Appendix B.    However,  Theorem~\ref{thm:Gram}  below implies that, unless one of these sets 
 is linearly dependent for all $t \in {\bf R} $, 
 this can happen for, at most, a finite number of values of $t$.
 
 \subsubsection{Generic in $t$ for $V_k$ fixed}
 
      \begin{thm}   \label{thm:Gram}
   Let  $V_1,  V_2, \ldots V_d $ be a fixed set of unitary  matrices  in $  M_d({\bf C}) $ and
   define $A_m$ as in \eqref{tform} with $t \in (-1,1)$.   Then either  $\{ A_m^* A_n \} $ is a linearly
   dependent set for all $t \in {\bf R}$ or it is a linearly independent set except for
  a finite number of  values  of $t \in (-1,1) $.   The same holds for  $\{ A_m A_n^* \} $.
     \end{thm} 
     \noindent{\bf Proof:}  We present two arguments.  
      
    (i)  With $A_m$ generated from $V_m$ as in \eqref{tform},  let $G(t) $ denote the Gram matrix   in 
    part (D) of Remark~\ref{rmk:prelim}.
   Then  $\det G(t) $  is a polynomial in $t$ and  $\{ A_m^* A_n \} $ is a linearly independent set 
     if  and only if $t$ is not a root of $\det G(t) = 0$.  The maximum degree  is determined
     by the diagonal which has $d$ elements  $g_{kk,kk} $  of order $t^4$  and  
     $d^2 - d$ elements $g_{jk,jk} $ ($j \neq k$ ) of order $t^2$. Thus $\det G(t)$  is a polynomial of degree at most
     $4d + 2(d^2 -d ) = 2d(d+1)$. 

 (ii) Every  $d \times d $ matrix $B$ with elements $b_{jk}$  can be associated with a
   $1 \times d^2 $ vector whose elements are 
   $v_{ d (j-1)  + k }  =  b_{jk} $.   Let $F(t) $ be the  $d^2 \times d^2 $ matrix whose
   rows are given by the vectors associated with the matrices $A_m^*  A_n $ in this way.
   Then   $\{ A_m^* A_n \}$ is a linearly independent set if and only if $\det F(t)  \neq 0 $.
    Now $ A_m^* A_m $ has exactly one element $t^2 $ (and all others $0, 1 $).   When $m \neq n$
   the elements of  $ A_m^* A_n  $  are of order $0$ or $1$ in $t$.
   Thus, $\det F(t) $ is a polynomial of degree  $\leq 2d + d^2 - d = d(d+1) $.  
   Since  $t = \pm 1$ are always roots, the maximum number of distinct roots
 in $(-1,1) $ is $ d^2 + d - 2 = (d+2)(d-1) $. \qed   
 
    Although the polynomials $ \det G(t) $ and $\det F(t) $ are not identical, they must have the same roots.
  However, these roots need not have the same degeneracy.  Numerical work suggests that the
  number of distinct roots   which lie in $(-1,1)$ is often much less than   $(d +2)(d-1) $ 
    because some roots are degenerate, imaginary, or lie outside  $(-1,1)$.
  For the example in Section~\ref{subsect:key},  there is a single highly degenerate root at
  $t = \tfrac{-1}{d-1} \in (-1,1)$  for each integer $d \geq 3$. For the example in Section~\ref{sect:odd},
  there are no roots in $(-1,1) $ when $d > 3. $   
  
  The simplest example of a situation in which $\det G(t) =  \det F(t) = 0 ~~\forall~ t $ is when all $V_m$ in \eqref{tform}
  are  diagonal unitaries.    At the end of Section~\ref{sect:odd}, we show that
   when $ d = 2 \nu $ is even and all $V_m $ are given by \eqref{eventry}, the $A_m$ generated from \eqref{tform}
  have $\det G(t) =  \det F(t) = 0 ~~\forall~ t $.
  Other cases are discussed in Section~\ref{sect:band}.
   
  \subsubsection{Algebraic geometry preliminaries}
 
  We now present some  results from algebraic geometry, which are needed for the rest of this section.
  All algebraic varieties and manifolds in what follows will be assumed to be real. 
  The terminology \emph{algebraic manifold} is used for an algebraic variety which is also a smooth manifold, i.e., without singularities.

A measure $\mu$ on an algebraic manifold $M \subset \mathbf{R}_m$ of dimension $d \leq m$ is said to be \emph{locally equivalent} to the $d$-dimensional Lebesgue measure if for each open subset $V$ of $M$ which is diffeomorphic to an open ball $B$ in $\mathbf{R}_d$, the restriction of $\mu$ to $V$ is equivalent to the
restriction of  Lebesgue  measure to $B$, i.e., they have the same null sets.

  \begin{lemma} \label{lm:alggeom} 
Let $M \subset \mathbf{R}_m$ be an algebraic manifold of dimension $d $, and let $N$ be an algebraic sub-variety of $M$ of dimension at most $d-1$. Let $\mu$ be a measure on $M$ which is locally equivalent to the Lebesgue measure. Then $\mu(N) = 0.$
\end{lemma}

\noindent{\bf Proof:} The set $N_1$ of  singular points of $N$ is an algebraic variety  of dimension at most $d-2$ (if not empty), and it is a   Zariski closed subset of $M$. For $x \in N \setminus N_1$, choose an open neighborhood $V_x \subseteq M$ of $x$ diffeomorphic to an open ball in ${ \bf R}_d$, and such that $V_x \cap N$ is diffeomorphic to a smooth submanifold in $V_x$ of dimension $d' = \dim(N)<d$. Since the $d$-dimensional Lebesgue measure of any smooth submanifold in ${ \bf R}_d$ of dimension less than $d$ is zero, and since the restriction of $\mu$ to $V_x$ is equivalent to the $d$-dimensional Lebesgue measure, we deduce that  $\mu(V_x \cap N)=0$. Since $N \setminus N_1 \subseteq \bigcup_{j=1}^\infty (V_{x_j} \cap N)$, for some countable set $\{x_j\}$ of points in $N \setminus N_1$ (because $N \setminus N_1$ is a Lindel\"of space, being second countable), we conclude that $\mu(N \setminus N_1)=0$.

Let $N_2$ be the set of singular points of $N_1$, and define, succesively, $N_{j+1}$ to be the set of singular points of $N_j$. Then $N_k$ is empty for some $k \le d+1$. Repeating the argument above, with $N_{j+1} \subset N_j \subset M$, $j=1,2, \dots, k-1$, in the place of $N_1 \subset N \subset M$, we see that $\mu(N) = \mu(N_1) = \mu(N_2) = \cdots = \mu(N_k) = 0$, as desired.   \qed
\medskip

    Let ${\cal U}(d)^m $ denote the set of $m$-tuples $(U_1, U_2, \ldots U_m ) $ of unitary matrices  $U_j \in  M_d({\bf C} )$.  
    \begin{lemma} \label{lm:null-set}
Let $d, k, m\ge 1$ be integers, and let $P \colon \mathcal{U}(d)^m \to \mathbf{C}$ be a function that arises from evaluating a polynomial in the $2md^2$ real variables given by the real and imaginary parts of the entries of the elements in $\mathcal{U}(d)^m$. If $P$ is not identically zero, then
$$Z := \big\{(U_1,U_2, \dots, U_m) \in \mathcal{U}(d)^m : P(U_1,U_2, \dots, U_m) = 0 \big\}$$
is a null set with respect to the normalized Haar measure on $\mathcal{U}(d)^m$. 
\end{lemma}
\noindent{\bf Proof}: 
Let  $(U_1,U_2, \dots, U_m)$ be an $m$-tuple of matrices with each $U_j \in M_d(\mathbf{C})$.  Consider the natural embedding
of this $m$-tuple into $\mathbf{R}_{2md^2}$, obtained by taking the real and the imaginary part of each entry of each $U_j \in M_d(\mathbf{C})$. 
Such an $m$-tuple  belongs to $\mathcal{U}(d)^m$ if it is in the zero set of finitely many polynomials in these $2md^2$ real variables. Hence  
$ \mathcal{U}(d)^m$ is a (real) algebraic variety. Moreover, it is an algebraic manifold (i.e., it has no singularities), because it is homogeneous, being a group. Thus, if it had one singularity, then all points would be singularities by homogeneity, which is impossible. Moreover, $\mathcal{U}(d)^m$ is connected in the usual Euclidean topology, and hence also in the Zariski topology. This implies that  $\mathcal{U}(d)^m$  is an irreducible algebraic manifold.

It is a standard fact from algebraic geometry that $Z$ is an algebraic sub-variety of $\mathcal{U}(d)^m$ of dimension strictly less than the dimension of
 $\mathcal{U}(d)^m$. It therefore follows from Lemma 3.1 that the Haar measure of $Z$ is zero. \qed

\subsubsection{Almost all channels are extreme}

We are now ready to prove the key result which implies that almost all channels generated using \eqref{tchan} and \eqref{tform} are extreme.
Although both (a) and (c) are special cases of (b) in the theorem below, we consider them
separately for ease of exposition.

\begin{thm} \label{thm:density}
Let $d \ge 3$ be an integer and fix $t \in (-1,1)$.   Suppose that there is a 
unitary matrix $ W  $ in $\mathcal{U}(d \! - \! 1)$ such that when all $V_m = W$ 
 in \eqref{tform}, then the resulting set  $\{A_m^*A_n\}_{n,m=1}^d$  is linearly independent.
 Then we can conclude the following:

\begin{itemize}
\item[\rm{(a)}] The set of $d$-tuples $(V_1, V_2, \dots, V_d)$ in $\mathcal{U}(d-1)^d$ for which the associated matrices $A_1,A_2, \dots, A_d \in M_d(\mathbf{C})$ defined in (10) satisfy that $\{A_m^*A_n\}_{n,m=1}^d$ is linearly independent is a co-null set in $\mathcal{U}(d-1)^d$ with respect to the Haar measure on $\mathcal{U}(d-1)^d$.

\item[\rm{(b)}] Fix a partition $\{J_1,J_2, \dots, J_\kappa \}$ of $\{1,2,\dots, d\}$. For each set of unitaries $W_1,W_2, \dots, W_\kappa $ in $\mathcal{U}(d-1)$, let $(V_1,V_2, \dots, V_d) \in \mathcal{U}(d-1)^d$ be given by $V_m = W_j$ when $m \in J_j$, and let $A_1, A_2, \dots, A_d$ be the matrices in $M_d(\mathbf{C})$ associated to $V_1,V_2, \dots, V_d$. The set of $\kappa$-tuples $(W_1,W_2, \dots, W_\kappa )$ in $\mathcal{U}(d-1)^\kappa $ for which $\{A_m^*A_n\}_{n,m=1}^d$ is linearly independent is a co-null set in $\mathcal{U}(d-1)^\kappa $ with respect to the Haar measure on $\mathcal{U}(d-1)^\kappa $.

\item[\rm{(c)}] For each unitary $V \in \mathcal{U}(d-1)$, let $A_1,A_2, \dots, A_d \in M_d(\mathbf{C})$ be as defined in  \eqref{tform} with all $V_m = V$.  
The set of $V \in \mathcal{U}(d-1)$ for which the set $\{A_m^*A_n\}_{n,m=1}^d$ is linearly independent is a co-null set in $\mathcal{U}(d-1)$ with respect to the Haar measure on $\mathcal{U}(d-1)$. 

\end{itemize}
 If, instead, $\{A_m A_n^* \}_{n,m=1}^d$  is linearly independent, the conclusions above hold for $\{A_m A_n^* \}_{n,m=1}^d$. 
 
\end{thm}

 \noindent{ \bf Proof:} (a) Let $G$ be the Gram matrix defined in part (D) of Remark~\ref{rmk:prelim}.
The function $P \colon \mathcal{U}(d-1)^d \to \mathbf{C}$ given by $P(V_1,V_2, \dots, V_d) = \det G$ is  a polynomial in the real and imaginary entries of the $V_m$, as in Lemma~\ref{lm:null-set}. A $d$-tuple $V_1, V_2, \dots, V_d$ belongs to the null-set $Z$ of $P$ if and only if the set $\{A_m^*A_n\}$ is \emph{not} linearly independent. By assumption, there exists a unitary $W$ such that $P(W,W, \dots, W) \neq 0$. Lemma~\ref{lm:null-set} therefore implies that $Z$ is a null-set with respect to the Haar measure.

(c) Let $G$ again be the Gram matrix defined in part (D) of Remark~\ref{rmk:prelim}. Then the  function $P \colon \mathcal{U}(d-1) \to \mathbf{C}$   
 given by   $P(V)  = \det G $ is a polynomial in the real and imaginary entries of  $V$, as in Lemma  \ref{lm:null-set}.  
 By assumption, there exists a unitary $W$ such that $P(W) \neq 0$. Lemma~\ref{lm:null-set} therefore implies that $Z$ is a null-set with respect to the Haar measure.
 
 (b)  Although (a) and (c) are special cases of (b), we proved them first to avoid cumbersome notation. The proof of (b) is  similar and 
further details are omitted.    \qed

 It is worth noting that, in Theorem~\ref{thm:density}  above, (a) does not imply (c) because   the set of
 $d$-tuples of the form $(V, V, \ldots V) \in \mathcal{U}(d-1)^d$ is a null set  with respect 
 to the Haar measure on $\mathcal{U}(d-1)^k$.   For the same reason, (a) does not imply (b).

     For the next result observe that when   $| x \ket  $  is a vector on the unit sphere in   
     $ {\bf C}_d $, then the matrix $ 2\proj{x} - \rmid_d $ is unitary and self-adjoint.   
     \begin{thm}   \label{thm:proj}
  Let $d \ge 3$ be an integer and fix $t \in (-1,1)$.  Suppose that there is a vector 
  $| w \ket  $  on the unit sphere in   $ {\bf C}_{d-1} $ 
such that when  all $ V_m   = 2 \proj{w} - \rmid_{d-1}  $
 in \eqref{tform}, then the resulting set  $\{A_m A_n \}_{n,m=1}^d$  is linearly independent.    Then the following hold:

{\rm a) }  For almost all   $| x \ket  $ on the unit sphere in  $ {\bf C}_{d-1} $, 
   when  $ V_m   = 2 \proj{x} - \rmid_{d-1} ~~ \forall~m$
 in \eqref{tform} the resulting set  $\{A_m A_n\}_{n,m=1}^d$  is linearly independent.
 
{\rm b) } The conclusions  {\rm  (a) -- (c)}  of Theorem~\ref{thm:density} hold for both   $\{ A_m^* A_n \}$  and  $ \{ A_m A_n^*\}  $.
     \end{thm} 
     
 \noindent {\bf  Proof:} View the unit sphere $S_{2d-1}$ of $\mathbf{C}_d$ as a real submanifold of $\mathbf{R}_{2d}$. As $S_{2d-1}$  is the zero-set of a polynomial (in $2d$ real variables), it is a real algebraic variety. The sphere is also without singularities, being homogeneous, so it is a real algebraic manifold. Hence we can apply Lemma 3.6. 

For each fixed $| x \ket  \in S_{2d-1}$, set $V_m = 2 |x \rangle \langle x | - I_{d-1}$ and let $A_m \in M_d(\mathbf{C})$ be given as in (10), for $m=1,2, \dots, d$.
 The $(j,k)$th entry of $V_m$  is $2x_j \overline{x}_k - \delta_{jk}$, which is a polynomial of degree two in the real variables 
 $\mathrm{Re}(x_j), \mathrm{Im}(x_j), \mathrm{Re}(x_k), \mathrm{Im}(x_k)$. Hence each entry of $A_m$ is a polynomial (of degree two) with respect to these variables. 

Consider the polynomial $Q(x) = \mathrm{det}(G)$,  for $ | x \ket  \in  S_{2d-1}$, where $G$ is the Gram matrix defined in part (C) of Remark 3.4. Then, as in the proof of Theorem 3.8, 
$\{A_m^*A_n\}$ is linearly independent if and only if $Q(x) \ne 0$. By assumption, there exists $| w \ket \in 
S_{2d-1}$ with the corrresponding $\{A_m^*A_n\}$ linearly independent which implies that $Q(w)$ is not identically zero.
 Hence, by Lemma 3.6, the zero-set of $Q$ is a Lebesgue null-set. (By
the comments above Lemma 3.6, the standard Lebesgue measure on the sphere $S_{2d-1}$ is locally equivalent to the Lebesgue measure on $\mathbf{R}_{2d-1}$.)  This proves part (a) of the theorem. Part (b) is an immediate consequence of Theorem 3.8  and the fact that $A_m  = A_m^*  $.  \qed
\begin{remark}
We could also consider $V_m = 2 \proj{x_m} - \rmid_{d-1} $ in \eqref{tform} using different vectors $|x_m \ket $ on the unit sphere in
${\bf C}_{d-1} $ and prove results analogous to (a) and (b) in Theorem~\ref{thm:density}.
\end{remark} 

Let $d \geq 3$,  fix $t \in (-1, 1) $, and assume that we can find a unit vector  $| x\ket   \in  
 {\bf C}_{d-1}  $  such that the map  $\Phi $ given by  \eqref{tchan} with all $V_m =   2\proj{x} - \rmid_{d-1}  $ in
\eqref{tform}  is extreme in either the set of UCP  or  CPT  maps.  Since $V_m$ is self-adjoint, it then follows from 
part (B) of Remark~\ref{rmk:prelim}
and the theorems above that almost every choice of $(V_1, V_2,  \ldots V_d$) in \eqref{tform} generates a UCPT
map $\Phi$ that is extreme in both the set of UCP maps and the set of CPT maps, for all the scenarios described
in Theorem~\ref{thm:density}.   

Thus, we are motivated to find unit vectors  $| x \ket  $ which generate extreme
UCP (or CPT) maps in this way.  This is done in the following sections.

 \begin{itemize}
 
 \item  In Section~\ref{subsect:key} we introduce a vector $| w \ket $ for which the hypothesis
  in Theorem~\ref{thm:proj}  holds for all $d \geq 3$  if $ t \neq \tfrac{-1}{d-1} $.
 
 \item   In  Section~\ref{subsect:rankone}  we describe evidence that the  hypothesis in Theorem~\ref{thm:proj}  holds  for all  $t \in (-1,1) $
  when $d = 3,4, 5, 6, 7 $,   and  conjecture that it holds for all $d \geq 3$.
             
 \item   In  Section~\ref{sect:odd} we  introduce a different type of unitary $W = W^* \in {\cal U}(d-1) $  for which  the hypothesis  in 
 Theorem~\ref{thm:density} holds for all $t  \in (-1,1)$  when $ d \geq 5 $ is odd.
 
 \end{itemize} 
 
 Although Theorem~\ref{thm:density} only requires $W$ to be unitary, in all of our examples $ W = W^*$.
 Similarly, in all of our examples  which satisfy the hypothesis of Theorem~\ref{thm:proj}, 
 the unit vector $| x \ket   \in   {\bf C}_{d-1}   $  is in ${\bf R}_{d-1} $.
 
 If, however, $W$ in   Theorem~\ref{thm:density}  is unitary, but not self-adjoint, then our conclusions are more restricted.
 Let $d \geq 3$ and fix $t \in (-1,1)$.
  \begin{itemize}
 
 \item[a) ] If \eqref{tform} and \eqref{tchan} generate  a map $\Phi$ which is extreme in the set of
 CPT maps,  then  almost every corresponding choice of $(V_1, V_2,  \ldots V_\kappa)$ in \eqref{tform} generates a  
map $\Phi$ that is extreme in   the set of CPT maps, for all the scenarios described
in Theorem~\ref{thm:density}.   

\item[ b)  ]If \eqref{tform} and \eqref{tchan} generate  a map $\Phi$ which is extreme in the set of
 UCP maps,  then  almost every corresponding choice of $(V_1, V_2,  \ldots V_\kappa ) $ in \eqref{tform} generates a  
map $\Phi$ that is extreme in   the set of UCP maps, for all the scenarios described
in Theorem~\ref{thm:density}.   

\end{itemize}

  \subsection{Critical Examples}  \label{sect:crit}
  
  \subsubsection{Key example}   \label{subsect:key}
Recall that  $ | \one_d \ket $  denotes the vector whose elements are all $d^{-1/2} $  
and define $W_d  \in M_d({\bf C}) $ as
\be   \label{Vdef}
    W_d = 2 \proj{\one_d }  -\rmid_d  .
\ee
It is easy to verify directly that $W_d$ is unitary and self-adjoint; in fact, its eigenvalues are 
$-1$ with multiplicity $d-1$, and $+1$ (non-degenerate).

For our first, and most important, example we choose $V_m$ in \eqref{tform} to be $W_{d-1}$ so
that     $A_1 = t  \proj{e_1} \oplus W_{d-1}$  and  
$ A_m =   S^{-1} A_{m-1} S =  S^{- m+1} A_1 S^{m-1}$ for $m = 2, 3, \ldots d $.  Then $A_m $ 
  can be written in block form as
\be  \label{keydef}
  A_m =    S^{-m+1} \pmx  t  & 0 \\ 0 & W_{d-1} \emx S^{m-1} , \qquad  m = 1, 2, \ldots d .
       \ee
For this example we will prove the following 
\begin{thm}  \label{thm:key}
For  $d \geq 3  $ and $t \in (-1,1) $,  let $A_m$ be the matrices defined in \eqref{keydef}.
Then  $A_m^* = A_m $ and when $  t \neq \tfrac{-1}{d-1} $,

{\rm a)}   the set of matrices $\{ A_m A_n \}_{m,n = 1}^d $ is linearly independent,  

{\rm b)}    the map $\Phi(\rho) =   \tfrac{1}{d-1 + t^2} \sum_{m = 1}^d  A_m \rho A_m $ is an extreme
 point of both  the CPT and UCP maps, and
 
{\rm c)}  $\Phi$ is not factorizable.
\end{thm}
Since $A_m = A_m^* $, part (B) of Remark~\ref{rmk:prelim} implies that part (b)  follows immediately from part (a). Part (c) then follows from (b)
because, as  stated in part (C) of Remark~\ref{rmk:prelim},  an extreme point of the UCP maps is never
factorizable.

The proof  of (a) is postponed to Section~\ref{sect:key},  in part because the argument is fairly long,  but also
because in the cases $d = 3$ and $t = \tfrac{-1}{d-1} $ we prove some related results of independent
interest.  The proof of   (a) is given in Sections \ref{sect:over4} -- \ref{sect:espace}
for $d \geq 4$ and in Section~\ref{sect:d=3} for $ d = 3$.   
For $d = 3$, we also prove that the channels
for $t = 1$ and  $t = \tfrac{-1}{d-1} = - \half $ have exact factorizations which are dual
in the sense that they can be obtained from the same unitary operator in $M_3({\bf C}) \ot  M_3({\bf C}) $
by switching the roles of the two algebras. 
 In Section~\ref{sect:d=3}, we consider the special case $t = \tfrac{-1}{d-1}  $ and show that
 the anti-symmetric and symmetric sets
$  \{ A_m A_n - A_n A_m \}_{m < n } $  and $  \{ A_m A_n + A_n A_m \}_{m < n } $ are each separately linearly dependent. 
 
\begin{remark} \label{rmk:more}
We also note that our proof (presented in Section~\ref{sect:key} ) implies the following

{\rm a) }
 For $ d \geq 3$ and  $t = 1$, the sets $\{ A_m^2 \}_{m = 1}^d $ and  
 $\{ A_m A_n  - A_n A_m\}_{m < n} $ are each separately linearly dependent.

{\rm b) }    For $ d \geq 3$ and   $t = - 1$, the set  $\{ A_m^2 \}_{m = 1}^d $ is linearly dependent, but
   the set $\{ A_m A_n \}_{m \neq n} $ is  linearly independent.   
   \end{remark}

 \subsubsection{Rank one projections} \label{subsect:rankone}
 
 The unitary $W_{d-1} = 2 \proj{ \one_{d-1} }  -  \rmid_{d-1} $ used in the previous section is invariant under 
 permutations, i.e., $P^{-1} W_{d-1} P = W_{d-1}$ for all permutations $P$.  This symmetry
allows us to present a full analysis of the linear independence of $\{ A_m A_n \} $ in Section~\ref{sect:key}.
However, this same symmetry gives rise to a highly degenerate linear dependency when $t = \tfrac{-1}{d-1} $.

  If $| \one_{d-1} \ket $ is  replaced by any unit vector $| x \ket \in {\bf C}_{d-1} $, then   $V = 2 \proj{x}   - \rmid_{d-1} $ 
 is also unitary. One would expect that most   unit vectors $ | x \ket $  give a unitary $V$
 which, when used in \eqref{tform},  generates a set of matrices for which  
  $\{ A_m A_n \} $ is  linearly independent   when $t =       \tfrac{-1}{d-1} $.  
 
\begin{conj}  \label{conj:new} 
For all integers $d \geq 3  $, one can find a vector $| x \ket $ on the unit sphere in $ {\bf R} _{d-1}  $  such that with 
$V_m =   2 \proj{x} - \rmid_{d-1}  ~~ \forall ~ m$  in \eqref{tform}, 
 the set of matrices $\{ A_m A_n \} $ is linearly independent when $t = \tfrac{-1}{d-1} $.
\end{conj}

 Exact numerical results obtained using Maple for  $d = 3, 4 ,5, 6, 7  $ not only support this conjecture, but suggest
 that it is easy to find such $| x \ket $.   Indeed, any vector for which $x_j \neq 0 ~~ \forall ~j $ 
and  some $ |x_j| \neq  (d-1)^{-1/2}$ seems to satisfy this conjecture.
 Moreover, when $d$ was small enough to find all roots of $\det G(t) = 0 $, there was
  high degeneracy at $t = \pm 1$, and many complex roots.  In view of Theorem~\ref{thm:Gram},
  it hardly seems plausible that every   $| x \ket  \in {\bf C}_d $ generates a Gram matrix from
  $\{A_m A_n \} $ which has a root at  $t = -\tfrac{-1}{d-1} $.
  
  % \pagebreak
 
  \subsubsection{$d = 2 \nu + 1  > 3 $ odd}   \label{sect:odd}

In this section we consider a completely different example 
which is not associated with a rank one projection.   Instead, it is based on a permutation matrix
with $\nu = \half(d-1)$ swaps, which can be viewed as a projection of rank $\nu = \half(d-1)$.
  
 All addition of indices in what follows will be mod $ d $.  For $m = 1, 2, \ldots d$, define
\be    \label{oddform}
          V_m  =   \sum_{k = 1}^\nu  \big(  |e_{m + k }  \kb e_{m- k} |   + |e_{m - k }  \kb e_{m + k} | \big) = 
         \Big(   \sum_{j=1 }^d   | e_j \kb e_{2m - j } |  \Big) -  \proj{e_m} .
         \ee
Then $V_m =V_m^* $ and its restriction to ${\bf C}_d \setminus \spn\{ | e_m \ket  \} \simeq {\bf C}_{d-1} $ is a unitary whose effect on a vector $| w \ket $
is to swap  its elements  in $\nu $ pairs  $ w_{j + \nu} \leftrightarrow w_{j - \nu } $, and  $A_m = V_m \op   t \proj{e_m} $ 
has the form \eqref{tform}.

 \begin{thm}   \label{thm:odd}
When $  A_m =   V_m  \op t \proj{e_m} $, then $A_m = A_m^* $, the set $\{ A_m^* A_n \}_{m,n = 1}^d $
 is linearly independent and 
the channel  $\Phi(\rho) = \sum_{k =1}^d A_k \rho A_k $  is an extreme point of both the
set of CPT maps and the set of UCP maps

a)   for all $t \in (-1,1)$  when $ d  = 2 \nu + 1 > 3 $, and 

b)  for all $t \neq -\half   \in (-1,1) $    when $d = 3$.
\end{thm}
\noindent{\bf Proof:} 
It follows from Proposition \ref{prop:diag} and the fact that   $A_m = A_m^*$   that it suffices
 to show that the set $\{ A_m A_n \}_{m \neq n } $ is linearly independent.  Observe that
 \be
    \big(V_m + \proj{e_m} \big) \big(V_n  +  \proj{e_n}\big)  & = &  \sum_{j = 1}^d   \sum_{k= 1}^d | e_j \kb e_{2m - j}| e_k\kb e_{2n - k } |    \nn \\
       &   =  &   \sum_{j = 1}^d  |e_j \kb e_{2(n-m) + j }|  = S^{2(n-m)}
       \ee
     where $S = \sum_k |e_k \kb e_{k+1} |$ is the cyclic shift.  It then follows that 
\begin{align}
    A_m A_n   &  =      S^{2 (n-m)}- (1-t) \big(  | e_{2m-n} \kb  e_n |    +  |e_m \kb e_{2n-m} | \big)  + \delta_{mn} (1-t)^2  \proj{e_m}   \nn  \\
\intertext{ or, equivalently, with $\ell = n - m $}  \label{oddcond}
 A_m A_{m+\ell}   &    =     S^{2 \ell }- (1-t) \big(  | e_{m -   \ell} \kb  e_{m + \ell} | + |e_m \kb e_{m + 2 \ell} | \big)  + \delta_{\ell 0}(1-t)^2  \proj{e_m}  .
\end{align}
  This implies that   $\{ A_m A_n \}_{m \neq n} $ can be decomposed into $d-1 $ disjoint sets of
     the form  ${\cal C}_{\ell} = \{ A_m A_{m + \ell } \}_{m = 1}^d $ with $\ell =  1, \ldots (d-1)$, i.e., 
     $\spn\{  {\cal C}_k \}  \cap \spn\{  {\cal C}_{\ell } \}  = \{ 0 \}  $ when $k \neq \ell  $.
     Thus it suffices to show that each of the  sets ${\cal C}_{\ell}  $  contains $d$ linearly independent matrices.

     For $\ell \neq 0$, the matrices in ${\cal C}_{\ell }$ are non-zero only where $S^{2\ell } $ is non-zero
     and     $\sum_{m= 1}^d A_m A_{m + \ell} = \big[ d-2(1-t) \big] S^{2 \ell } $    (which is non-zero unless $d = 3, t = -\half$).   Thus, 
      it suffices to show that for each fixed $\ell \in \{ 1, 2,\ldots  d-1 \} $,  the matrices
     $$ F_{m \ell} \equiv    \tfrac{1}{1-t} \big(S^{2 \ell}  - A_m A_{m + \ell} \big)   =   | e_{m -   \ell} \kb  e_{m + \ell} | +  |e_m \kb e_{m + 2 \ell} | , \qquad  m = 1, 2, \ldots d$$
    are linearly independent.   Now,
     after a shift which does not affect linear independence, we
     can associate each   $F_{m \ell}  $ with the vector   $| v_{m \ell } \ket = | e_ m \ket + | e_{m - \ell} \ket  $  in    ${\bf C}_d$
     where $m = 1, 2, \ldots d$ and  $\ell \in \{ 1, 2,\ldots  d-1 \} $ is fixed.
          
       One way to show linear independence is to consider  the matrix $M_\ell$ in which the vectors  $ |v_{m \ell} \ket  $ are the rows.
     Then $M_\ell = I + S^{-\ell} $  is  a circulant matrix whose eigenvalues are well-known to be
     $\lambda_k = 1 + e^{-2 \pi i k \ell/d } \neq 0 $   \cite[Section 2.5, Problem 21]{HJ} .    
     For $d$ odd, no $d$-th root of unity can be $-1$  which implies   $\lambda_k \neq 0$ and   $M_\ell$ is non-singular.
     Therefore, the set of
     matrices  $\{ F_{m \ell} \} _{m = 1 }^d$ or, equivalently, $\{ A_m A_{m + \ell } \} _{m = 1 }^d$ is linear independent for
    each fixed $\ell$. 
     
        Alternatively, observe that since $| e_j \ket = | e_{j - d \ell} \ket$,  one can write
    \bee
       2 | e_j  \ket = | e_j   \ket  + |  e_{j - d \ell}  \ket  & = &  \sum_{k = 1}^d (-1)^{k+1} |  e_{j - (k-1)\ell}  \ket  + |  e_{j - k \ell}  \ket   \\
          & = &   \sum_{k = 1}^d   (-1)^{k+1} |  v_{ (j - (k-1) \ell ) \, ,\, \ell} \ket .
    \eee
    This implies that, for each fixed $\ell$, every  vector in the orthonormal basis 
    $\{ | e_j \ket \}_{m=1}^d$ is in the span of $\{ | v_{m \ell }\ket  \}_{m =1}^d $ which implies that the vectors  $\{ | v_{m \ell }\ket  \}_{m =1}^d $
   form a basis for ${\bf C}_d$  and are,  hence,  linearly independent.   \qed   
 %  \medskip
   
   To see why this example is the opposite of a rank one projection, 
   recall that $E = E^* = E^2 \in {\bf C}_d$ is a projection of rank $r$ if and only if $ \rmid_d - E $ is a projection of rank $d - r $.
   Whenever $E$ is a projection, $V  = 2 E - I $ is unitary and  the replacement $E \mapsto  (I - E ) $  takes $V \mapsto -V $.
   Thus any unitary  matrix that can be constructed with a projection of  rank $> \lceil d/2 \rceil$ can also be constructed with a projection of  
   rank $\leq  \lceil   d/2 \rceil $.
   
  When $V_m $ is given by \eqref{oddform}, there is a projection $E_m = \sum_{k = 1}^\nu \proj{x_{mk} }$ of rank $\nu = \half(d-1)  $
   formed from the  vectors
    $|x_{mk}  \ket =  \1rt2 \big( |e_{m + k}  \ket  +  |e_{m-k  } \ket  \big)$ so that  $V_m = 2 E_m - \rmid_{d-1}.  $
      Thus $V_m $ comes from a projection with the maximal    rank on   $ {\bf C}_{d-1}$ .

   When $d = 2 \nu $ is even, $d -1 $ is odd and it is not possible to construct a unitary
   consisting only of $\nu $ swaps.   If one uses \eqref{oddform}, the term with $k = \nu $
   becomes  $|m + \nu \kb m - \nu | =  |m + \nu \kb m +  \nu | $  ~ $(\bmod ~d = 2 \nu) $.
   When $d$ is odd,   $V_m $ is essentially a projection $W  \in {\bf C}_{d-1}  $  with $1$'s on the skew diagonal and $0$ elsewhere.
  Extending this to $d = 2 \nu $, i.e., generating $A_m $ in \eqref{tform} from   a projection $W  \in {\bf C}_{d-1}  $   with $1$'s on the skew diagonal, gives
   \be   \label{eventry} 
    V_m  =   \sum_{k = 1}^{\nu - 1}  \big(  |e_{m + k }  \kb e_{m- k} |  |e_{m + k }  \kb e_{m- k} | \big) + \proj{m + \nu },
    \ee
  However, when $A_m =  t \proj{e_m}  \oplus V_m $, 
   $$  S^{-\nu } A_m S^\nu   - A_m =  (1 - t) \big( \proj{e_m} - \proj{e_\nu } \big) $$
   which implies that  $A_m A_{m + \nu } $ is diagonal so that  
   $A_m A_{m + \nu }  \in  \spn \{ A_m^2 \} $.   There does not seem to be a natural
   generalization of \eqref{oddform} to $d = 2 \nu $   which yields a linearly independent set of $\{A_m A_n \} $.

\subsubsection{Band Width}  \label{sect:band}

We have  focused on finding examples of matrices $V \in {\bf C}_{d-1} $
that generate linearly independent sets of $\{ A_m^* A_n\}   $ via \eqref{tform}
because these imply that almost all choices of $V_k $ also generate  linearly independent sets
as described in Theorem~\ref{thm:density}.
Nevertheless, a null set with respect to Haar measure does allow for infinitely many choices of $V_k $
 that  generate linearly dependent sets.  We now consider  what those sets might look like.
 
 First, if all of the $V_m$ in \eqref{tform} are diagonal, then $A_m^* A_n $ is also diagonal.
 It then follows immediately from Theorem~\ref{prop:diag} and Remark~\ref{rmk:diag} that 
 $\{ A_m^* A_n \} _{m=1}^d  $  are linearly dependent.   In fact, if $A_m^* A_n $ is diagonal for even
 one pair of $m \neq n $, the set  of $\{ A_m A_n\}   $  is linearly dependent.  One can have
 such pairs even when none of the $V_m$ are diagonal.   This is the case for the
  $V_m $ in \eqref{eventry}, which
 generate  $A_m $ for which $A_m^* A_{m + \nu} $ is diagonal when
 $d = 2 \nu $ is even.
 
  We can extend the diagonal examples by introducing the  notion of {\em band width}.
  We first observe that when $S $ is the cyclic shift defined above \eqref{zeroform},
  every matrix  $B \in M_d({\bf C}) $ can be written uniquely as
  \be
         B = \sum_{ k = - \xi_-}^ {\xi_+}   D_k S^k
  \ee
 where $D_k$ is diagonal, $\xi_- =  \lfloor (d-1)/2 \rfloor$ and $\xi_+ = \lceil (d-1)/2 \rceil$.
 \begin{defn}  \label{def:band} 
 Let $d \ge 3$ be an integer. For each matrix $B $ in $M_{d}(\mathbf{C})$ with elements $b_{jk} $ define its  {\rm{\bf band width} } to be 
 $\beta(B) = \max \{  |j-k| : b_{jk} \ne 0\}$ (with $0 $ the  band width of the zero matrix). 
  Define its {\rm{\bf cyclic band width}}  $ \mu$ as the smallest positive integer such that   $B = \sum_{ k = -\mu}^{\mu}   D_k S^k $
  with each $D_k$ diagonal.
 \end{defn} 
 The following properties are straightforward to verify for all $A, B \in M_d({\bf C}) $.
 \begin{enumerate}[label=\alph{*})]
   \item  $\mu(B) \leq \beta(B) $     and    $\mu(B) \leq  \lceil (d-1)/2 \rceil$.
  
    \item  A matrix $D$ is diagonal if and only if $\mu(D)=0$ if and only if $\beta(D) = 0$.
    
    \item  $\beta(A^*) = \beta(A) $;   $\mu(A^*) = \mu(A) $ and  $\mu(S^*AS) = \mu(A)$; 
       
    \item  $\beta(AB) \le \beta(A) + \beta(B)$, and     
 $\mu(AB) \le \mu(A) + \mu(B)$.
    
 \end{enumerate}   
    We also note that one can find $B$ with $\beta(B) = d-1$.  Moreover, $\beta(B) $ is not invariant under
    cyclic permutations. Indeed, $S ^*|e_{d-1} \kb e_d | S =  |e_d \kb e_1 |$ so that a single cyclic shift
    can map the  matrix $B = |e_{d-1} \kb e_d | $ with $\beta(B) = 1$ to one with $\beta(S^* B S) = d-1 $, which is the maximal value of $\beta$.
    \begin{prop} \label{prop:band} 
    Let $d \ge 3$ be an integer and let $V_1, \dots, V_d$ be unitaries in $M_{d-1}(\mathbf{C})$ with band width $\beta(V_m) < (d-1)/4$, for each $m$. Then, for the associated matrices $A_1, \dots, A_d$ in $M_{d}(\mathbf{C})$, defined in (10), it follows that the set $\{A_m^*A_n\}$ is \emph{not} linearly independent for any value of $t \in [-1,1]$.
\end{prop}
\noindent  {\bf Proof:} Let 
$B_m =  t \proj{e_1}  \op V_m = \begin{pmatrix} t & 0 \\ 0 & V_m \end{pmatrix}, ~~ m=1,2, \dots, d.$
Then $\beta(B_m) = \beta(V_m) < (d-1)/4$  and $A_m = S^{-m+1} B_m S^{m-1}$  , which implies  $\mu(A_m) =\mu(B_m) < (d-1)/4$, for all $m$, so that
 $\mu(A_m^*A_n) < (d-1)/2$, for all $m,n$. Since the set of matrices $A$ in $M_d(\mathbf{C})$ satisfying $\mu(A) < (d-1)/2$ is a proper linear subspace of $M_d(\mathbf{C})$, we conclude that 
$\{A_m^*A_n\}$ cannot be a basis for $M_d(\mathbf{C})$, and hence not linearly independent.  \qed
  
   Note, however, that  $\mu(V_m)  >    (d-1)/4$ does  not necessarily imply that $\{A_m^*A_n\}$ generated as in \eqref{tform}  
   are linearly independent.  When $d  = 2 \nu \geq 8$,  the    $V_m $ in 
  \eqref{eventry}  have the maximal cyclic band width of $\nu = \lceil \half(d-1) \rceil   = \half d$, but
  nevertheless generate sets  of $\{A_m^*A_n\}$  which are linearly dependent for all $t \in {\bf R} $.

\section{Analysis of key example }   \label{sect:key}

\subsection{Overview}  \label{sect:over4}

We begin with an outline of the steps in the rather long analysis of the key example
in Section~\ref{subsect:key}.  In the next section we describe $A_m A_n$ explicitly
and find some properties which allow us  to reduce the linear independence of
$\{A_m A_n \}_{m \neq n } $ to that of $\{ X_{mn} \}_{m \neq n} $ where the 
matrices  $ X_{mn} $ each have only two non-zero rows and columns.

In Section~\ref{sect:eval} we observe that it suffices to consider the linear independence
of the sets  $\{ X_{mn }^+ \}_{m < n } $ and $\{ X_{mn }^- \}_{m < n } $   where
$X_{mn}^\pm = X_{mn} \pm X_{mn}^* $.   We then observe that we can remove a 
common factor so that the non-zero elements  $x_{jk} $ of $X_{mn }^\pm$ are  $\pm 1$
with the exception of $x_{mn} $ and $x_{nm} $.  The linear independence of the sets 
$\{ X_{mn }^\pm \}_{m < n } $ depends only on the elements above the diagonal.
We construct a larger matrix $\Omega^\pm $ whose rows are these elements
$x_{jk} $ of $X_{mn} $ arranged in lexicographic order.  Remarkably, the  
diagonal of $\Omega^\pm $ is a multiple of the identity,  whose value depends on $t$.
This allows us to reduce the linear
independence of  $\{ X_{mn }^\pm \}_{m < n } $ to an eigenvalue problem for $\Omega^\pm $.

The eigenspaces of $\Omega^\pm $ are associated with representations of the symmetric group,
which allows us to find all of the eigenvalues and eigenspaces explicitly.  This is done in
Section~\ref{sect:spec} for $d \geq  4$.  The results imply linear independence of $\{A_m A_n \}_{m \neq n } $ 
when $t \neq \tfrac{-1}{d-1} $.  

The reduction process used in Section~\ref{sect:structure} is not valid when $d= 3$ or when $t = \tfrac{-1}{d-1} $.
 Therefore the case $d = 3$ is analyzed separately in  Section~\ref{sect:d=3} in which we also present
a pair of novel  factorizability results relating the  channels with $t = 1$ and $t = -\half $.
The linear dependence of  $\{A_m A_n \}_{m \neq n } $ when $d > 3 $ and $t = \tfrac{-1}{d-1}$  is analyzed
in Section~\ref{sect:spec}.

\subsection{Structure of $A_m A_n $}   \label{sect:structure}

 We begin our analysis of our key example in Section~\ref{subsect:key}  by observing  that
 the entries of the matrix  $A_1 = t \proj{e_1} \oplus  W_{d-1} = t \proj{e_1} \oplus (2 \proj{ \one_{d-1} } - \rmid_{d-1} )$   are
\be   \label{key}
  a_{jk} = \begin{cases}  t  &  j = k = 1  \\
  v_{jk}  = \tfrac{2}{d-1} &  j \neq k \in  \{ 2, 3 \ldots d  \} \\
   0  &  j = 1 ~ \hbox{or} ~ k = 1,  j \neq k \\
                                     v_{jj} = \tfrac{3-d}{d-1}  & j = k \neq 1 \\
                                         \end{cases}  
                                        \ee 
Since $A_m = A_m^* $, $A_m A_n^* = A_m^* A_n = A_m A_n $ and a straightforward calculation  gives
\be     \label{AmAnels}
   (d-1)^2 \bra e_j , A_m A_n e_k \ket = \begin{cases} 
     - t (d-1)(d-3) & j = k = m,   j = k = n  \\
       2 t (d-1) & j = m, k \neq m,n  ~ \hbox{or} ~ k = n, j \neq m,n \\
      2 (d-3) & j = n, k  \neq m,n ~ \hbox{or} ~ k = m, j \neq m,n \\
          0  &  j = m, k = n \\
          4(d-2) & j = n, k = m  \\
         (d-3)(d+1)  & j = k \neq m,n  \\   
          -4 &  j \neq k,  ~~ j, k \neq  m, n  \end{cases}        
\ee
so that, e.g., 
  \bee    
   A_1 A_d  & = &  \frac{1}{(d-1)^2}   \pmx  \tau & b & \ldots& b  & \ldots &b & 0  \\ a &   & & &  &   &b \\  \vdots &&&&&& \vdots  \\
                a &    &   &  \wtd{V}_{d-2}^2  & &    &  b \\  \vdots &&&&&& \vdots  \\
                 a &   &&  & &  & b \\ u & a &   \ldots& a & \ldots  &  & \tau \emx
 \eee
 where $a = 2(d-3), ~ b = 2t(d-1),~  \tau =   - t (d-1)(d-3),~ u = 4(d-2) $  and 
 $  \wtd{V}_{d-2} $ is the matrix in $M_{d-2}({\bf C}) $ obtained by removing 
 the last row and column from $V$.
 
By Remark~\ref{rmk:diag}  we do not need to consider the diagonal elements of   
 $ A_m A_n $; therefore, it suffices to 
determine whether or not there exists  a matrix  $C$ with elements $c_{jk} $ such that   
      $\bra e_j ,   \sum_{m \neq n} c_{mn} A_m A_n   e_k \ket =  0 $ when $j \neq k $. 
       It  follows immediately from \eqref{AmAnels} that this holds   if and only if
          \be  \label{fullcond}
     0 = 4(d - 2) c_{kj} +   a \sum_{m \neq j, k}  (c_{jm} + c_{mk} ) + b \sum_{m \neq j, k} (c_{mj} + c_{km} )  - 4 \sum_{m,n \neq j,k} c_{mn}.
  \ee

We can simplify the conditions on $c_{jk} $ by first observing that it also follows immediately from
 \eqref{AmAnels} that % for all $j \neq k $
 \begin{align*}    \label{sum}
     (d-1)^2   \sum_{m \neq n}   \bra e_j , A_m A_n e_k \ket   & =  
               4 (d-2) [1 + t (d-1) ] &   \forall ~ j \neq k,   \\  
          \intertext{and}  
     (d-1)^2   \sum_{m \neq n}   \bra e_k , A_m A_n e_k \ket   & =   
           (d-1)(d-3) \big[ (d-2)(d+1) - 2 t (d-1) \big]  & \forall ~k
 \end{align*} 
  so that
 \be   \label{sumall}
    (d-1)^2   \sum_{m=1}^d \sum_{n =1}^d   A_m A_n =   p_d(t) \proj{  \one_d} + [q_d(t) -  d \cdot p_d(t)  ] \rmid_d
 \ee
 where
\ssq   \begin{align}
     p_d(t) & =   4 d (d-2) [1 + t (d-1) ]  \\
     \intertext{is independent of $j , k $ and non-zero unless $t = \tfrac{-1}{d-1}$, and}
     q_d(t) & =     (d-1)(d-3) \big[ (d-2)(d+1) - 2 t (d-1) \big] .
 \end{align}   \esq
 In the special case,
 $ t = \tfrac{-1}{d-1}$,
\be   \label{special}
  \sum_{m \neq n} A_m A_n  = \begin{cases}  0  & d = 3 \\  q_d \Big( \tfrac{-1}{d-1} \Big) \, \rmid_d  & d \geq 4 \end{cases} 
  \ee 
 with $q_d \Big( \tfrac{-1}{d-1} \Big)  =  d(d-1)^2(d-2)  $.   For $d = 3$,  
 this immediately implies that  $\{ A_m A_n  \}_{m \neq n } $ is linearly dependent when $ t = \tfrac{-1}{d-1} $

  If $t \neq \tfrac{-1}{d-1}$, then after using the freedom from Remark~\ref{prop:diag} to adjust the diagonal elements
arbitrarily, we can proceed as if   $\sum_{m  n } A_m A_n $ is a multiple of  $ \proj{ \one_d } $.
Then we can remove 
 $  \wtd{V}_{d-2}^2 $ by replacing $A_m A_n $ by  $(d-1)^2 A_m A_n + 4 d \proj{\one_d }$.
 Thus we can conclude that
 \begin{prop}   \label{prop:equiv}
For $d \geq 4$ with  $t \in   (-1,1) $ and $t \neq \tfrac{-1}{d-1} $, 
 the set $ \{ A_m A_n \}_{m, n = 1}^d $ is linearly independent  if and only if
  $\{ X_{mn} \}_{m \neq n}$ is linearly independent where
\be   \label{replace}
  X_{mn} = (d-1)^2 A_m A_n + 4 d \proj{\one_d } - D_{mn},   \qquad m, n = 1 , 2, \ldots d
    \ee  
    and $D_{mn}  = D_{nm}   $ is a diagonal matrix chosen so that the diagonal of $X_{mn}$ is identically zero.
    In particular, one can find a constant  $\lambda$ such that  
    \be    \label{Dmn}
   D_{mn}    = \lambda \rmid_d + (\tau - \lambda) \big( \proj{e_m} + \proj{e_n} \big) .
    \ee
    \end{prop}
  
    In Section~\ref{sect:spec}, we consider the case $t = \tfrac{-1}{d-1} $ in detail for $d \geq 4$.  Although   we can not write  $\proj{\one_d} $ as
    a linear combination of  $A_m A_n$ in that case, we can still use $X_{mn}$ to
    reach some conclusions about the linear independence of  $\{ A_m A_n \}_{m \neq n} $. Moreover, we give a simple proof of
    the linear dependence conditions in this case which does not use  Proposition~\ref{prop:equiv}.
  The case $d = 3$ must be treated separately, which is done in Section~\ref{sect:d=3}.

The matrix $X_{mn}$ obtained in \eqref{replace} has elements    
\bee
    \bra e_j ,  X_{mn}  e_k \ket = \begin{cases}  
       \wh{a} \equiv   2 (d-1) & j = n, k \neq m,n ~ \hbox{or} ~ k = m, j \neq m,n \\
       \wh{b}  \equiv    4 + 2t(d-1) &      j = m, k \neq m,n  ~ \hbox{or} ~ k = n, j \neq m,n \\
       \wh{u}  \equiv   4(d-1) &  j = n, k = m \\
       4 & j = m, k = n \\
       0 & \rm{otherwise} \end{cases} 
  \eee
or, equivalently,
\be    \label{xdef}
   X_{mn}  =    \wh{u}  |e_n \kb e_m |  +  4 |e_m \kb e_n | +  \wh{a}  \! \sum_{j \neq m, n} \!  \big( |e_j \kb e_m | + |e_n \kb e_j | \big) 
                     +   \wh{b} \!  \sum_{j \neq m, n} \! \big( |e_m \kb e_j | + |e_j \kb e_n | \big)      \quad             
\ee
 where   $\wh{a} = 2(d-1) ,  ~ \wh{u} =  4(d-1)$,   ~ $\wh{b}  = 2t(d-1) + 4 $.  
Thus we can write,  showing only  rows and columns with non-zero elements,   
\bee
X_{nm} ~   =  \quad \begin{blockarray}{c  c c c c c c c c c}
 & & & & & m & & n & & \\  \\
\begin{block}{cc (cccccccc)}
&  &   &  &  & \wh{a} &   & \wh{b} &  & \\
& &   &  &  & \vdots &  & \vdots & &  \\
& &   &  &  & \wh{a} &     & \wh{b} &  &  \\
m &  & \wh{b} & \ldots & \wh{b} &0 & \wh{b} & 4 & \wh{b} & \ldots \\
& &  & & & \wh{a} &   & \wh{b} & & \\
& &  &  &  & \vdots &   & \vdots &  &  \\
& &  &  & \ & \wh{a} &   & \wh{b} & \ &  \\
n &  & \wh{a} & \ldots & \wh{a} & \wh{u} & \wh{a} & 0 & \wh{a} & \ldots \\
& &   & &  & \wh{a} &  & \wh{b} &  &  \\
& &   &  &  & \vdots &  & \vdots &  &  \\
\end{block}
\end{blockarray}   ~.
\eee
It follows from \eqref{xdef} that the set $\{ X_{mn} \}_{m \neq n}  $ is linearly dependent if
and only if there is a matrix $C$  with elements $ c_{mn} $ such that
\be      \label{lindep}
   0 & = & \sum_{mn} c_{mn} \bra e_j , X_{mn} e_k \ket   \qquad \forall ~~  j, k  .
   \ee
\begin{remark}
The set of matrices satisfying \eqref{lindep}  or, equivalently, \eqref{fullcond} is a subspace ${\cal N}$ of $M_d({\bf C})$ which is invariant under 
permutations and under the adjoint map, i.e.,
  $C \in {\cal N}  $ implies $C^* \in {\cal N}$   and  
$   P^* C P \in {\cal N}  $  for all permutation matrices   $P$.

\end{remark} 
We chose ${\cal N}$ to denote this subspace because it is the null space
of the $d^2 \times d^2 $  matrix whose rows are the elements of $X_{mn} $.
However, this formulation does not give much insight.  Instead we will exploit
 invariance  under the adjoint map to reformulate the linear dependence as pair of eigenvalue problems.
Moreover, the invariance under permutations will allow us to identify the eigenspaces
with irreducible representations of the symmetric group $S_d$.

\subsection{Reformulation as an  eigenvalue problem}  \label{sect:eval}

We can replace  the pair
  $  X_{mn} , X_{nm} $ with $m < n$ by the pair of matrices 
   $$X_{mn}^{\pm} \equiv X_{mn } \pm X_{mn}^ * = X_{mn} \pm X_{nm} $$
  so that  $(X_{mn}^+)^* = X_{mn}$ is symmetric and  $(X_{mn}^-)^*  = -  X_{mn}^- $ 
  skew symmetric.  Since
  $ \tr X_{mn}^+ X_{jk}^-  = \tr (X_{mn}^+ X_{jk}^- )^T = -  \tr X_{mn}^+ X_{jk}^-  =  0 ~~ \forall  ~ j < k, $ and $ m < n$,
    the sets $ \{ X_{mn}^+ \}_{m < n}  $ and   $\{ X_{mn}^- \}_{m < n} $ are
  orthogonal.  Thus the set  $\{ X_{mn} \}_{m \neq n } $ is linearly independent if and only if both
of the sets $ \{ X_{mn}^+ \}_{m < n}  $ and   $\{ X_{mn}^- \}_{m < n} $ are  linearly independent.

  For $d \geq 4$, we observe that when $t = \frac{d-3}{d-1 } $,  $  \wh{a} = \wh{b} $ and
   $$X_{mn}^- = - ( \wh{u} - 4) \big( |e_m \kb e_n |  -  |e_n \kb e_m | \big) = 4(2-d)  \big( |e_m \kb e_n |  -  |e_n \kb e_m | \big)$$
  from which the linear independence of  $\{ X_{mn}^- \} $ follows immediately.   For all other choices of $t$,
  $  \wh{a} - \wh{b} \neq 0 $ and   $  \wh{a} + \wh{b} \neq 0 $ for any choice of $t$.  Thus, in what follows, we
  can assume that   $  \wh{a} \pm  \wh{b} \neq 0  $, and
modify the definition of $X_{mn}^\pm $ above by removing the common
non-zero factor of  $ (  \wh{a} \pm \wh{b} ) $  and defining  (with $t \neq  \frac{d-3}{d-1 } $ for $X_{mn}^-$)
\be
 X_{mn}^\pm  & \equiv &  \frac{1} {\wh{b} \pm \wh{a} }\big(   X_{mn } \pm X_{mn}^ * \big) =   \frac{1} {\wh{b} \pm \wh{a} } \big( X_{mn} \pm X_{nm}  \big) \\ \nn ~~ \\ \nn
  & = &    w_d^\pm \big(  |e_m \kb e_n |  \pm  |e_n \kb e_m |  \big)  
  +  \sum_{j \neq m, n}   \big( |e_j \kb e_m | + |e_n \kb e_j |      \pm   \big( |e_m \kb e_k | \pm  |e_k \kb e_n | \big)  
\ee
where we used \eqref{xdef} and 
\bee
    w_d^+(t)  & = &   \frac{  \wh{u} + 4}{ \wh{a} + \wh{b} }     =      \frac{2d}{d + 1 +t (d-1) }  \qquad   ~~ \in  (1,d)   \qquad  \qquad -1 < t < 1 \\  
   w_d^-(t)    & =   &   \frac{  \wh{u} -4 }{ \wh{a} - \wh{b} }    =  \frac{2(d-2)}{   (d-3) - t(d-1)  }  \qquad  
    \in   \begin{cases} (1 , \infty) & -1  < t < \tfrac{d-1}{d-3} \\ (-\infty, 2-d)  &  \frac{d-3}{d-1}  <  t < 1 \end{cases} 
 \eee
  Note that
 $w_d^\pm $  depends on $t$ and when this is important we write $w_d^\pm(t) $;
 otherwise we suppress this dependence to avoid cumbersome notation.  In
 particular
 \begin{align*}
   w_d^+(-1) & = d  &    w_d^+ \big( \tfrac{-1}{d-1} \big) & = 2  &   w_d^+(0) & = \tfrac{2d}{d+1}  &  w_d^+ \big( \tfrac{d-3}{d-1} \big) & = \frac{d}{d-1} & w_d^+(1)  & = 1  &\\
   w_d^-(-1) & = 1 &    w_d^- \big( \tfrac{-1}{d-1} \big) & = 2  &   w_d^-(0) & = 2 \tfrac{d-2}{d-3}  & w_d^- \big( \tfrac{d-3}{d-1} \big) & = \pm \infty & w_d^-(1)  & = 2 -d  &  
 \end{align*}
Although the replacement \eqref{replace}   does not preserve linear independence when $t = \tfrac{-1}{d-1}$,
we allow $w_d^\pm  = 2$ in  order to obtain a complete set of eigenvectors for
 the eigenvalue problem for $\Omega_d^\pm(0) $ described below.
  
  It will be useful to consider $X_{mn}^\pm $ as a function of $x$, where $x$ replaces $w_d^\pm(t) $,
  and allow $x \in {\bf R} $ even when there is no $t \in[-1,1]$ for which $x = w_d^\pm(t) $.  Thus,
          \be      X_{mn}^\pm(x)  
           =  x  \big( |e_m \kb e_n |   \pm    |e_n \kb e_m | \big)     
      +  \!  \sum_{j \neq m , n } \! \!  \big( \, | e_m \kb e_j |  \pm   | e_j \kb e_m| +  |e_j \kb e_n | 
         \pm     |e_n \kb e_j | \, \big)  
    \ee
or, equivalently, again showing only the non-zero rows and columns. 
\bee
X_{nm}^\pm (x) ~   =  \quad \begin{blockarray}{c  c c c c c c c c c}
 & & & & & m & & n & & \\  \\
\begin{block}{cc (cccccccc)}
&  &   &  &  & \pm 1&   & 1&  & \\
& &   &  &  & \vdots &  & \vdots & &  \\
& &   &  &  & \pm 1&     & 1&  &  \\
m &  & 1& \ldots & 1&0 & 1& x & 1& \ldots \\
& &  & & & \pm 1&   & 1& & \\
& &  &  &  & \vdots &   & \vdots &  &  \\
& &  &  & \ & \pm 1&   & 1& \ &  \\
n &  & \pm 1& \ldots & \pm 1& \pm x & \pm 1& 0 & \pm 1& \ldots \\
& &   & &  & \pm 1&  & 1&  &  \\
& &   &  &  & \vdots &  & \vdots &  &  \\
\end{block}
\end{blockarray}
\eee
 
 We now define a $ \half d (d-1 ) \times  \half d (d-1 )$ matrix  
  $\Omega_d^\pm(x) $ whose rows are given by the elements of $ X_{mn}^\pm(x) $
 above the diagonal, arranged in lexicographic order so that its elements are
 as follows with $m < n ,  j < k  $
\be      \label{vmndef}
  \omega_{mn,jk}  & = &  ( x \mp 2) \delta_{jm} \delta_{kn} + \delta_{jm} + \delta_{mk} \pm  \delta_{jn} \pm  \delta_{nk}        \nn \\
          & = & \begin{cases}   x  & j = m, k = n \\  +1 &  j = m, k \neq m,n  \hbox{~or~ }  k = n, j \neq m,n \\
            -1 & j = n, k \neq m,n  \hbox{~or~ }  k = m, j \neq m,n \\   0 & \hbox{otherwise} \end{cases}
  \ee
 Since $\Omega_d^\pm(x) =  \Omega_d^\pm(0)  + x \rmid_d $   we can conclude that
  \begin{thm}    \label{thm:eval}
The elements of each of the sets $\{ X_{mn}^\pm(w_d^\pm) \}_{m < n} $ are linearly independent if and only if  $ - w_d^\pm$
is not an eigenvalue of $\Omega_d^\pm(0) $.
 \end{thm}
 
 We will show that
  \begin{prop} \label{thm:negeval}
 The eigenvalues of  $\Omega_d^-(0) $ are
    
    $ \bullet ~~$  $d-2$ with multiplicity $d-1$, and 
    
    $\bullet ~~$  $-2$ with multiplicity $\half (d-2)(d-1) $. 
 \end{prop}

  \begin{prop}    \label{thm:poseval}
   The eigenvalues of $\Omega_d^+(0) $ are 
   
     $ \bullet ~~ $  $2(d-2)$ non-degenerate, 
     
     $ \bullet ~~ $  $d-4 $ with multiplicity $d-1$, and 
    
     $ \bullet~~  $  $-2$ with multiplicity $\binom{d-1}{2}  - 1 = \half d(d-3) $. 
     
      \end{prop}
   
   These results were conjectured after using Mathematica for $d = 3, 4, 5, 6$.  We  prove
   them by giving,   for each eigenvalue, a  linearly independent set of eigenvectors
   that is at least as large as
   the claimed multiplicity.  Since the total number of eigenvectors  of $\Omega_d^{\pm} $ can not be larger
than  its dimension, which is $\half d(d-1) $,  both lists of eigenvalues above are exhaustive.   
    
\subsection{Description of the eigenspaces for $d \geq 4$}  \label{sect:espace}

Since each eigenspace of $\Omega_d^\pm(0) $ is the null space of $\Omega_d^\pm(x) $ for
the eigenvalue $-x$, we will describe them as these null spaces.   
  Although each eigenvector is an element  of ${\bf C}_{ d(d-1) /2} $, we  write them as the
  associated $d \times d$
   symmetric or skew symmetric matrices (both with zero diagonal), in the case of 
   $\Omega_d^+ $ and $\Omega_d^-$ respectively.  This facilitates associating each eigenspace with an
   irreducible representation of  $S_d $. 
  It also allows  us to translate  linear dependence relations  for elements of the sets $\{ X_{mn}^\pm \}_{m < n} $
   to relations for elements of $ \{ A_m \pm A_n \}_{m < n} $, which we will write using commutators and anti-commutators
   defined, respectively,  as
   \be
       [ A, B] = AB - BA   \qquad  \qquad  \{ A, B \} = AB + AB  .
   \ee
   
   For each $x$ for which $\Omega_d^\pm(x) $ is singular, we will exhibit enough linearly
independent eigenvectors of the null space to demonstrate that the multiplicity of the
eigenvalue $-x$ is at least as large as claimed above.  Although
 Proposition~\ref{prop:equiv} also does not hold when $t =  \tfrac{-1}{d -1} $ (which corresponds to $x = 2$)
 the associated eigenspaces are still needed to complete the proofs of Propositions~\ref{thm:negeval} and \ref{thm:poseval}.
 Moreover, as  explained in Section~\ref{sect:spec}, the eigenspaces we find can still
   be associated with linear dependence of certain subsets of $\{ A_m A_n \} $.   

 Because Proposition~\ref{prop:equiv} does not hold  for $d = 3$,  that case is analyzed separately  in Section~\ref{sect:d=3}
(where it is also observed that the results below do hold.)  Therefore, in what follows we assume that $d \geq 4$.

 {\bf Skew symmetric, $x = 2-d =  w_d^-(1)$:}  
The null space is spanned 
by skew-symmetric matrices $C_k$, $k = 1, 2,  \ldots d$ with  
\be
         C_k = \sum_{j \neq k} \big( |e_{k} \kb e_j | - |e_j \kb  e_k | \big ) .
\ee
One readily verifies that $\sum_{k = 2}^{d} C_k = - C_1 $ and that $C_k$ with $ k = 2, 3, \ldots d $
give a basis of $d-1$ linearly independent vectors associated with the $d-1$ Young tableaux
 $\begin{array}{c}~ \\ \young(1~~~~~,k) \end{array}$.  This corresponds to the
 linear dependence 
 \be   \label{skew:t=1}
     \sum_{m \neq n }  [ A_m A_n ]  = 0   \qquad   \hbox{for each} ~   n = 1 ,2 \ldots d
 \ee
 when  $t = 1$, and proves part (a) of Remark~\ref{rmk:more}.

 {\bf Skew symmetric, $x = 2 = w_d^-\Big( \tfrac{-1}{d-1} \Big) $:}    
 Although  Proposition~\ref{prop:equiv} does not hold in this case,
  $X_{mn}^-  $ is a multiple of
$[A_m , A_n]$ so that the  null space of $\Omega_d^-(2) $ 
 describes the linear dependence of the matrices  $\{ A_m A_n - A_n A_m \}_{m < n} .$  
     In this case,
the null space is spanned by $\half (d-1)(d-2) $ skew symmetric $d \times d $ matrices 
  $C_{jk}$  with $ 1 < j < k \leq d$  given by
 \be   \label{C2espace}
 C_{jk} =  |e_1\kb e_j | - |e_j \kb e_1| -   |e_1 \kb e_k| +  |e_k \kb e_1|  + |e_j \kb e_k | - |e_k \kb e_j | .
 \ee 
 Since the term  $ |e_j \kb e_k | $ occurs in one and only one $C_{jk} $ with  $ 1 < j < k \leq d$  
 these $\half  (d-1)(d-2)$ matrices are clearly linearly independent.  In fact, they are associated with
Young tableaux of the form
 $\begin{array}{c}~ \\ \young(1~~~~~,j,k) \end{array}$.    For  $t = \tfrac{-1}{d-1}$  this generates 
the  linear dependence relations 
 \be   \label{skewx=2}
         [A_j, A_k ] + [A_k , A_m ]    + [ A_m, A_j ]  = 0  
 \ee
with  $j, k, \ell $ distinct.

 \medskip
 
  {\bf Symmetric, $x = 2(2-d) $: }  The null space of $\Omega_d^+\big(2(2-d)\big) $ is readily seen to be  
    $\proj{\one_d} $ (or $d \proj{\one_d} - \rmid_d $ if one wants to keep the diagonal zero) which 
 corresponds to the trivial
 representation of $S_d$  with Young tableaux $ \young(12~~~d)  $.
 This gives
    $ \sum_{m \neq n}  X_{mn}^+ = 0 $ when $x = 2 (2-d) $.  Although $w_d^+(t) \neq 2(2-d)$ when $t \in [-1,1]$,
   this reflects the fact that $\sum_{m \neq n} A_m A_n $ has the form \eqref{sumall}.
    \medskip
 
{\bf Symmetric, $x = 4-d $: }  $w_d^+(t) = 4-d$ in the domain $[-1,1]$ under consideration only in the case $ d= 3 $,
which is studied in Section~\ref{sect:d=3}.   For $d = 4$,  $w_4^+(t) \neq 0$  for any 
 $t \in {\bf R}$ and for  $d \geq 5$,    $w_d^+(t) = 4-d $ only for $t \in [-4,-1) $.
  The null space of $\Omega_d^+(4-d) $ is spanned by the symmetric matrices
\bee
  C_{k \ell } =  \sum_{j \neq k, \ell } \big( |e_k \kb e_j | + |e_j \kb  e_k | - |e_\ell \kb e_j | -  |e_\ell \kb e_j |  \big)     \qquad  k \neq \ell .
  \eee
The matrices $C_{1 k }$ with $k = 2, 3, \ldots d$ are readily verified to be linearly independent,
 giving  a basis of $d-1$  matrices
corresponding to the $d-1$ Young tableaux
 $\begin{array}{c}~ \\ \young(1~~~~~,k) \end{array}$.     For $d \geq 4$,    
     this null space corresponds  to linear dependence relations  
     \be
          \sum_m  \big( X_{jm}^+ - X_{jm}^+  \big) = 0 
     \ee
for each fixed choice of $j < k $.  Then \eqref{replace} implies
 \be    \label{4-sym}
 \sum_{m \neq j, k }   \big(  \{ A_j , A_m \} - \{ A_k , A_m \}  \big) = \wtd{D}_{jk} 
   \ee
where $ \bra e_n  , \wtd{D}_{jk}  \, e_n \ket = 0 $ if $n \neq j,k $.  \medskip

{\bf Symmetric, $x = 2 =  w_d^+\Big( \tfrac{-1}{d-1} \Big) $:}    
  The null space of $\Omega_d^+(2) $  is spanned by 
symmetric matrices of the form $C_{jk,mn}^+ = B_{jk,mn} + B_{jk,mn}^* $, where $j, k, m ,n $
are distinct  in $\{ 1, 2 \ldots d\}$, and 
\be
  B_{jk,mn} \equiv |e_m \kb e_j | -   |e_m \kb e_k |  - |e_n \kb e_j | +   |e_n \kb e_k | 
\ee
and can be identified with the  irreducible representation of $S_d$ described by the
Young diagram  $\begin{array}{c}~ \\ \yng(7,2) \end{array}$  which can readily be shown
(using standard hook length arguments, e.g.,
\cite[Section 2.8]{S}) to have dimension $\half d(d-3) $.
We can verify this independently by observing that the matrices
\be     \label{basis2sym} 
     \{ B_{2k,1n} :  3 \leq  n < k \leq d \}  \cup  \{ B_{2k,13} :  k = 4, 5 \ldots d \}
\ee
are linearly independent.  There are $\sum_{k = 4}^d (k-3) = \half (d-3)(d-2) $ matrices
in the first group and $d-3$ in the second for a total of $\half d(d-3) $  matrices.
(Or one can observe that each $k$ is associated with a total of $k -2$ matrices
and $ \sum_{k = 4}^d (k-2)  = \sum_{j = 2}^{d-2} j = \half (d-2)(d-1) - 1 $.)
The matrices in \eqref{basis2sym} do not belong to the null space; however, the set of adjoints of the matrices
in \eqref{basis2sym} also form a linearly independent set whose elements are 
linearly independent of those in \eqref{basis2sym}.  Hence the corresponding sets
$ \{ C_{j,mn}^\pm \equiv B_{jk,mn} \pm B_{jk,mn}^*  \}$  are  each linearly independent.  Thus we
have found that the set
\bee
      \{ C_{2k,1n}^+ :  3 \leq n < k \leq d \}  \cup  \{ C_{3k,12}^+:   k = 4, 5 \ldots d \}
\eee
 gives  $\half d(d-3) $ linearly independent elements of the null space of $\Omega_d^+(2) $.
 Note that this basis is formed from the Young tableaux  
 $\begin{array}{c}~ \\ \young(12~\ldots~,nk) \end{array}$ and 
  $\begin{array}{c}~ \\ \young(13~\ldots~,2k) \end{array}$   in standard form.
  This gives linear dependence relations of the form
  \be   \label{sym:X}
       X_{jk}^+ - X_{k,m}^+ + X_{m,n}^+ - X_{n,j}^+  = 0
  \ee
which hold for any choice of $j,k,m,n$ distinct.

% \bigskip

%  \pagebreak

\subsection{Analysis for $d = 3$}  \label{sect:d=3}

When $d = 3$, the results of the previous section hold formally, however, 
the arguments given there are not valid because they depend on Proposition~\ref{prop:equiv}
which does not hold for $d = 3$.    Moreover,  in addition to  showing that a channel with
$t \in (-1,1) $ is extreme unless
 $t = \tfrac{-1}{d-1} = - \tfrac{1}{2}$, we show that when $t = - \half $ it is factorizable.
Furthermore, the channel with $t = 1$ has an exact  factorization through $M_3({\bf C})  \ot M_3({\bf C}) $ which 
uses a different unitary $\bU$ than in  \eqref{factpm1}, and
which is dual to the  factorization for $t = - \half $  in the sense 
 that these channels can be obtained  using the same unitary conjugation
 after exchanging the roles of the subalgebras.

When $d =3$, instead of  \eqref{replace}, we simply define 
$X_{mn}^\pm =  \pm \tfrac{ 1}{t} (A_n A_m \pm A_m A_n ) $, 
in which case $  \pm \tfrac{1}{t} $ plays the role of $w_3^\pm(t)$.
 Then  $  \Omega_3^\pm(x)  = \pmx  x  & 1 & \pm 1 \\ 1 & x  & 1 \\  \pm 1 & 1 & x  \emx  $
so that   
$\det \Omega_3^\pm(x)  = x^3  - 3x \pm 2 = (x \pm 2)(x \mp 1)^2 $,
which is consistent with  Propositions~\ref{thm:negeval}  and \ref{thm:poseval} above.
In particular,

\begin{itemize}
       
                 \item  For $t = - \tfrac{1}{2}, x = -2$, the matrix $\Omega_3^+(-2) $ has a one-dimensional null space.

               \item  For $ t = - \tfrac{1}{2},  x = + 2$, the matrix $\Omega_3^-(+2) $ has a one-dimensional  null space.   
                                          
        \item For  $ t = 1, x =   1 = 4 - d $, the matrix  $\Omega_3^+(+1) $ has a two-dimensional  null space.     
           
          \item  For $ t  = 1, x =   - 1 = 2-d $,  the matrix $ \Omega_3^-(-1) $ has two-dimensional  null space.            
\end{itemize} 

We first consider the case $t = \tfrac{-1}{d-1} =  - \half $,  for which the situation for $d = 3$ differs slightly from that  for $d \geq 4$.
\begin{prop}   \label{prop:d=3-half}
       For $ d = 3 $ and $t =  \tfrac{-1}{d-1} = -\half $,  the set   $\{ A_m^2 \}_{m = 1}^d $ is linearly
       independent, but
   
   a)   $\sum_{m \neq n} A_m A_n =   \sum_{m < n }\{ A_m,  A_n    \} = 0 $,  which implies that 
   $  \{ A_m A_n + A_n A_m \}_{m < n } $ is linearly  dependent,   
   
   b)  $  [A_1,  A_2]  +  [A_2 , A_3 ] + [ A_3 , A_1 ]  = 0$,  which implies that
     $\{ A_m A_n - A_n A_m \}_{m < n } $ is  linearly dependent,  and
   
   c) $A_1 A_2 + A_2 A_3 + A_3 A_1  =  A_1 A_3 + A_3 A_2 + A_2 A_1 = 0$.   

\noindent Moreover,  the map $\Phi(\rho) =  \tfrac{4}{9} \ds{ \sum_{m = 1}^3  A_m^* \rho A_m }$ has an exact factorization through  $M_3({\bf C}) \ot M_3({\bf C}) $.
\end{prop} 
\noindent{\bf Proof:}   Part  (a) follows immediately from \eqref{special}.  One can also observe that the
null space of  $ \Omega_3^+(-2)$  is spanned by $|e_1 \kb e_2 | -  |e_1 \kb e_3 | + |e_2 \kb e_3 |  $
  which implies (a).  
To prove (b) observe that  the null space of  $ \Omega_3^-(+2) $ is spanned by
$|e_1 \kb e_2| - |e_2 \kb e_3 | + |e_3 \kb e_1 | $.
  Part  (c) follows immediately from (a) and (b).   The final assertion follows immediately
  from   Theorem~\ref{thm:duald3}  below.  \qed

To distinguish the cases $t = -\half $ and $t = 1$, we now  use
$B_k$ to denote the matrices associated with $t = 1$.
The matrices $B_k$ are unitary which implies that $\{ B_m^2 \}_{m = 1}^3 $ are linearly
dependent.    However, as noted above, $\Omega_3^\pm(x) $ has a pair 
of two dimensional null spaces when $t = 1$, 
which leads to additional linear dependence relations.     Thus, we find
 \begin{prop}  \label{prop:d=3t=1}
  For   $ d = 3$ and  $t = 1$, and   $j ,k, \ell$ distinct
      
     {\rm a)}     $ \{ B_j, B_k  \} - \{ B_j  B_\ell \}  = 0 $, which implies that 
   $  \{ B_m B_n + B_n B_m \}_{m < n } $ is linearly  dependent,  
       
        {\rm b)}   $     [B_j, B_k ]   + [ B_j , B_\ell ] = 0   $,   which implies that
     $\{ B_m B_n - B_n B_m \}_{m < n } $ is  linearly dependent,  and
         
{\rm c)}      $   B_1 B_2    =  B_2 B_3 = B_3 B_1 =   Q$  and  $B_2 B_1 = B_3 B_2 = B_1 B_3 = Q^T $ 
  with   $Q =  \tfrac{1}{4} \pmx   0 & 0 & 1 \\ 1  & 0 &0  \\ 0 & 1 & 0 \emx $.   
  
 \noindent Moreover, the channel  $ \Psi(\rho) =   \tfrac{1}{3}  \ds{ \sum_{j=1}^3  B_j \rho B_j }$ has an exact
  factorization  through $M_3({\bf C}) \ot M_3({\bf C} ) $ which uses a different unitary  $\bU$ than in
   Remark~\ref{rmk:factpm1}.
   
  \end{prop} 
  \noindent{\bf Proof:}   To show (a),  it suffices to  observe that  the null space of $\Omega_3^+(1) $ is spanned
  by the pair of  vectors
  $|e_1 \kb e_2 | - |e_1 \kb e_3 | $ and  $|e_1 \kb e_2 | - |e_2 \kb e_3 | $.
  To show (b), it suffices to  observe that  the null space of $\Omega_3^-(-1) $ is spanned
 by the   pair of   vectors $| e_1 \kb e_2 | +  | e_1 \kb e_3 | $
   and $| e_1 \kb e_3 | +  | e_2 \kb e_3 | $.
     The equivalences in part (c) follow immediately from parts (a) and (b),    
and simple computation  or  \eqref{AmAnels} gives $Q$.   
(Note that  (a) and (b)   are exactly what one would get by formally using
\eqref{skew:t=1} and \eqref{4-sym} when $d = 3$.)  The final assertion follows immediately from the next result.

\begin{thm}   \label{thm:duald3}
  The channels $\Phi$ and $\Psi$, defined in part (d) of    Propositions~\ref{prop:d=3-half} and \ref{prop:d=3t=1}
respectively, have exact  factorizations through $M_3({\bf C} ) \ot M_3{(\bf C} )  $  which are dual in the sense that there
is a unitary matrix   $\bU  \in M_3({\bf C} ) \ot M_3{\bf C} ) $
such that 
\begin{align}    \label{UfactPhi}
  \Phi(\rho)    = \tfrac{4}{9} \sum_{k = 1}^3  A_k^* \rho A_k  & =  ({\cal I}_3 \ot  {\rm Tr} ) \big( {\bf U}^* (\rho \ot \tfrac{1}{3} \rmid_3  ){\bf U}   \big), \quad  and  \\
    \label{UfactPsi}
  \Psi(\rho)   =   \tfrac{1}{3}   \sum_{k=1}^3  B_k^* \rho B_k  & = ( {\rm Tr}  \ot {\cal I}_3   ) \big( {\bf U}^* (\rho \ot \tfrac{1}{3} \rmid_3  ){\bf U} \big)~.
  \end{align} 
\end{thm} 
\noindent{\bf Proof:}   We first define
\begin{align}   \label{Udef} 
     {\bf U}  &  =  \tfrac{2}{3} \pmx A_1  & A_3   & A_2 \\ A_3 & A_2 & A_1 \\ A_2 & A_1 & A_3 \emx   \qquad  = 
               \tfrac{2}{3} \sum_{k=1}^3 A_k \ot B_k  \\
   \label{Wdef}    { \bf W}   &  =   \tfrac{1}{3}  \pmx  - B_1 & 2 B_3 & 2 B_2  \\  2 B_3  & - B_2  & 2 B_1 \\
                  2 B_2 & 2 B_1 & - B_3 \emx       =   \tfrac{2}{3}  \sum_{k = 1}^3 B_k \ot A_k   ~.
\end{align}
 It follows immediately from part (c) of  Propositions~\ref{prop:d=3-half} and \ref{prop:d=3t=1}
respectively, 
that both  $ {\bf U}$ and ${\bf W} $ are unitary, and that   \eqref{UfactPhi} and \eqref{UfactPsi} hold.  \qed

We note that since  $ {\bf U}$ and ${\bf W} $ are identical except for the exchange  $A_k \leftrightarrow B_k $
we also have
\begin{align*}   %\label{WfactPhi}
 \Phi(\rho)     =      ( \tr \ot {\cal I}_3   ) \big( {\bf W}^* (\rho \ot \tfrac{1}{3} \rmid_3  ){\bf W} \big)    \qquad   \hbox{and} \qquad 
      \Psi(\rho)    =       ({\cal I}_3 \ot {\rm Tr} ) \big(  {\bf W}^* (\rho \ot \tfrac{1}{3} \rmid_3  ){\bf W}  \big) .  
\end{align*} 
The interplay between the channels for $t = -\half $ and $t = 1$ is interesting.  They are dual in the
sense that they can be obtained    by simply switching the subspace  over which 
one takes the partial trace of $ {\bf U}^* (\rho \ot \tfrac{1}{3} \rmid  ){\bf U} $ with the  same unitary 
for both channels.

\subsection{The special case   $t = \tfrac{-1}{d-1}$ with $d \geq 4$.}   \label{sect:spec}

\subsubsection{Linear dependence relations}

In this case, $p_d(t) = 0 $ so that \eqref{sumall} does not imply that $\proj{\one_d}$ is a linear combination of 
elements of $\{ A_m A_n  \}$.  Therefore, it is not clear that we can draw any conclusions about the linear independence
of the set  $\{ A_m  \pm A_n  \}_{m < n} $ from that of $\{ X_{mn}^\pm \}_{m < n} $.   However,
the resulting linear dependence relations for $\{ A_m  \pm A_n  \}_{m < n} $ can also be proved directly,
as observed after the following
\begin{prop}  \label{prop:d=4}
For $ d \geq 4$ and  $t =  \tfrac{-1}{d-1}$, the set   $\{ A_m^2 \}_{m = 1}^d $ is linearly independent, but
  $\sum_{m \neq n } A_m A_n =   4  \rmid_d =  \tfrac{9}{7}   \sum_{m=1}^d A_m^2 $, which implies
that  $ \{ A_m A_n  \}_{m,n = 1}^d $ is linearly dependent.  Moreover, each of the sets $\{ A_m A_n + A_m A_n \}_{m < n} $
and $\{ A_m A_n - A_m A_n \} _{m < n } $ is linearly dependent, and the following relations hold  
for any choice of $j, k, m, n$ distinct
 \begin{subequations}    \label{jkmn} \be     
  \label{skew:jkl}     0   & = &      [A_j, A_k ] + [A_k , A_m ]    + [ A_m, A_j ]   \\
 \label{skew:jkmn}    0 & = &       [A_j, A_k ] + [A_k , A_m ]    + [ A_m, A_n ]  + [A_n , A_j ]  \\
  \label{sym:jkmn}    0 & = &       \{A_j, A_k \} - \{A_k , A_m \}    + \{ A_m, A_n \}  - \{A_n , A_j \} \\
  \label{spec:jkmn}   0 & = &   (A_j - A_m) ( A_k - A_n ) .
\ee   \end{subequations} 
\end{prop} 
\noindent{\bf Proof:}  It is straightforward to verify \eqref{spec:jkmn} directly.   For example,   observe that for $d = 5$,
the matrices $  A_1 - A_2 $ and $A_5 - A_4 $ are, respectively, 
\bee
    \tfrac{1}{4} \!  \pmx     d \! - \!  4 & 0 & - 2 & -2  &  -2  \  \\  0 & 4 \! - \!  d & 2 & 2 &   2 \\
                  -2 & 2 & 0   & 0  &   0 \\  -2 & 2 & 0   & 0  &   0 \\
                  -2 & 2 & 0   & 0  &   0    \emx       \qquad \qquad 
                   \tfrac{1}{4} \! \pmx     0   & 0  &   0 &   2 & - 2  \\  0   & 0  &   0 &  2 & - 2   \\
                 0   & 0  &   0 &   2 & - 2  \\   2 & 2 & 2 & 4  \! - \!  d & 0  \\  - 2 & - 2& - 2 & 0 &  d  \! - \! 4 \emx .
 \eee
Then
$ (A_j - A_m) ( A_k - A_m )  \pm   ( A_k - A_m )(A_j - A_m) = 0 $,   gives  
 \eqref{sym:jkmn} and \eqref{skew:jkmn}, respectively.   Moreover,
 \eqref{skew:jkl} can be obtained from    linear combinations of permutations of
\eqref{skew:jkmn}.   For example, let  $X_{jkmn} $  denote the expression on the right
in \eqref{skew:jkmn}.  Then  $X_{1234} + X_{1324} + X_{1243} $ yields \eqref{skew:jkl}
with $j,k,m \simeq 1,2,4$.

Although the argument above does not use the results from Section~\ref{sect:espace}, it is
worth describing the connection.  
In the skew symmetric case,      \eqref{replace} implies that   
 $  [ A_m , A_n]  =  \tfrac{\wh{b} - \wh{a}}{(d-1)^2} X_{mn}^-$ so that   \eqref{skewx=2} implies
\eqref{skew:jkl}.   Then, by using \eqref{skew:jkl} with  $j,k,m $ and with
 $j,m,n $ together with the fact that $[A_j,A_m] = - [A_m, A_j ] $, one can prove  \eqref{skew:jkmn}. 
  
The symmetric case requires slightly more work.    Because the signs in \eqref{sym:X} alternate, it follows 
immediately  from  \eqref{replace} that the terms with $\proj{\one_d}$ cancel so that
\bee     
   \{ A_j, A_k  \} - \{ A_k, A_m \} + \{ A_m , A_n \} - \{ A_n, A_j \} =  D_{jk,mn} 
\eee 
where  $    D_{jk,mn}  $ is a diagonal matrix.   Moreover, it follows from \eqref{Dmn} that
 $    D_{jk,mn}  $  is identically zero. 
Combining this with \eqref{skew:jkmn} implies
$   (A_j - A_k)(A_m - A_n ) =  D_{jk,mn}  = 0  $.  (Instead of using \eqref{Dmn}, one could  
  use  (49d)   to conclude that  $ D_{jk,mn} = 0$.)  

Note that, as described  after \eqref{C2espace} 
and \eqref{basis2sym} respectively,  there
are $\half(d-1)(d-2) $   independent constraints of the form \eqref{skew:jkl}  or \eqref{skew:jkmn}
and    $\half d(d-3) $ independent 
 constraints of the form   \eqref{sym:jkmn}.  In addition, 
   \eqref{special}  implies that $\sum_{m \neq n } A_m A_n $ is a multiple of $\sum_{j=1}^d  A_j^2 $.
 Thus, the subspace of $M_d({\bf C} )$ spanned by $\{ A_m A_n \}_{m,n =1}^d $ has dimension
 $  3d -2 $.  
 
%When $d = 4$, there are 3 independent constraints of the form  \eqref{skew:jkl}  or \eqref{skew:jkmn} and two 
% of the form \eqref{sym:jkmn}.  A direct calculation 
%  gives   $   \{ A_j, A_k  \}  + \{ A_m , A_n \} = \tfrac{4}{3} \rmid_4 $ for $j,k,m,n$ distinct, 
%  which implies \eqref{sym:X}.  Moreover, $ \bU =   \tfrac{3}{4}  \pmx  A_j - A_k & A_m - A_n \\  A_m - A_n & A_j - A_k \emx$ 
%  is unitary because  \eqref{spec:jkmn}  implies that the columns are orthogonal  and  
%    $\sum_\ell A_\ell ^2  - \{A_j,A_k\} -   \{ A_m,-A_n \}  = \tfrac{28}{9} \rmid_4 - \tfrac{4}{3} \rmid_4 $.
%    However, this does not imply factorizability because    $(\id_4 \ot \tr) \bU (\rho \ot \half \rmid_2  ) \bU  \neq  \Phi(\rho) $.     
             
             \subsubsection{Factorizability when $d = 4$.}                                              

 The very high level of linear dependence when $t = \tfrac{-1}{d-1} $ makes it natural to conjecture that the associated channels 
 are factorizable, as is the case when $d = 3$.  Therefore, we considered this question when $d = 4$.  
 Although we did not fully resolve it, we have some partial results.  Moreover, we can show that
a necessary condition is that $Y_j Y_j^*$  and $Y_j^* Y_j $ in Propositon \eqref{prop:factcond} are 
  both independent of $j$.  This suggests that   {\em if} the channel  is factorizable, then $Y_j = U_j$ is unitary.
 In this section, we give the conditions on $U_j$ and discuss their implications.
 
 \begin{remark}  \label{rmk:d4fact}
For  $d = 4, t = - \third  $ and  $A_k$  given by \eqref{keydef}, let $\Phi$ be the
 channel  in part~(b) of Theorem~\ref{thm:key}.  
Then $\Phi$ has an exact factorization through $M_d(\mathbf{C}) \otimes \mathcal{N}$ with 
$\bU = \sum_{k= 1}^4 A_k \otimes U_k $ with $U_k \in {\cal N} $ unitary and $\tau(U_j U_k^* ) = \delta_{jk} $
if and only if the following  conditions hold:
\begin{align}   \label{jk=mn}
    Q_{jk}^+ =  Q_{mn}^+    & \quad   \hbox{ and }  \quad  R_{jk}^+ =  R_{mn}^+  &   j,k,m,n  ~~\hbox {distinct} \\
    \label{U+cond}
       Q_{jk}^+ + Q_{km}^+  + Q_{mj}^+ = 0    &  \quad  \hbox{ and }  \quad    R_{jk}^+ + R_{km}^+  + R_{mj}^+ = 0 &   j,k,m  ~~\hbox {distinct} \\
\label{U-cond} 
       Q_{jk}^- + Q_{km}^-  + Q_{mj}^- = 0    &  \quad  \hbox{ and }  \quad    R_{jk}^- + R_{km}^-  + R_{mj}^- = 0 &  j,k,m  ~~\hbox {distinct} 
 \end{align}
with  $Q_{jk}^\pm \equiv U_j U_k^* \pm U_k U_j^* $ and $R_{jk}^\pm   \equiv U_j^* U_k \pm U_k^* U_j $. 
 \end{remark}

Before giving the proof, we observe  that \eqref{U+cond} and \eqref{U-cond}   are equivalent to the asymmetric pair of conditions
\be  \label{asym:jkl}
       U_j U_k^* + U_k U_m^* + U_m U_j^* = 0     \qquad  and  \qquad    U_j ^*U_k + U_k^* U_m + U_m^* U_j  = 0.
\ee
Adding   \eqref{U+cond} and \eqref{U-cond} gives  \eqref{asym:jkl}.  Conversely, by combining \eqref{asym:jkl}
with its adjoints, we can  recover \eqref{U+cond} and \eqref{U-cond}.
  
\noindent{\bf Proof:}  We begin without the assumption that $Y_j $ is unitary and  (with a slight abuse of notation) let 
$Q_{jk}^\pm \equiv  Y_j Y_k^* \pm  Y_k Y_j^*  $.
When $s \neq t $, the conditions in \eqref{factcond}  for $Y_j Y_k^* $  are
formally the same as those on $c_{jk} $ in \eqref{fullcond}.  Thus, a necessary condition for factorizability is  % (with  $j,k,m,n $ distinct)
\be  \label{zz1a}
     4 Y_k Y_j^*   =  2 Q_{mn}^+  +  \sum_{\ell \neq j,k} Q_{k \ell} ^-- Q_{j \ell} ^-   \qquad \qquad     j ,k,m,n  \hbox{~~distinct}.
       \ee
Then $Q_{j \ell}^- = - Q_{\ell j}^-$  implies that $4  (Y_k Y_j^* + Y_j Y_k^* ) = 4 Q_{jk}^+ = 4 Q_{mn}^+  $ as in \eqref{jk=mn}.
 Next, observe that $ Y_j Y_k^* - Y_k Y_j^* $ gives  $ T_{jk} \equiv 2 Q^-_{jk} + Q^-_{km} + Q^-_{mj} - Q^-_{jn} - Q^-_{nk} = 0 $.
Then taking  $T_{jk} + T_{mn} $ (after  $j \mapsto m,  k \mapsto n, m \mapsto k, n \mapsto j$) we obtain
\be  \label{U-jkmn}
   Q^-_{jk} + Q^-_{km}  + Q^-_{mn} + Q^-_{nj}  = 0   \qquad  \qquad  j,k,m,n  ~~\hbox {distinct} 
\ee
which has the same form as \eqref{skew:jkmn}.  Then, just as \eqref{skew:jkmn} implies \eqref{skew:jkl}, \eqref{U-jkmn}
formally implies \eqref{U-cond}.

 When $s = t$, the conditions  in \eqref{factcond} become (for each fixed $j = 1,2,3,4$)
\be   \label{s=ta}
    28  I_\mathcal{N}  & = & Y_j Y_j^* + 9 \sum_{k \neq j } \Big(  Y_k Y_k^* +    Q_{jk}^+ +  5   Q_{mn}^+ \Big)   \qquad  \quad m,n \hbox{~ distinct from } j,k  \nn  \\
         & = & Y_j Y_j^* + 9 \sum_{k \neq j }  \Big( Y_k Y_k^* +  6   Q_{jk}^+      \Big)                   \ee 
 Taking the difference between \eqref{s=ta} with $j= 1$ and  with $j = 2$ implies
\be
    0 =  8   (Y_2 Y_2^*-  Y_1 Y_1^*)  +   6 ( Q_{13}^+  + Q_{14}^+ - Q_{23}^+  - Q_{24}^+ )  =   8   (Y_2 Y_2^*  -  Y_1 Y_1^*).
\ee
Repeating for other pairs of $j \neq k $, we conclude that $Y_j Y_j^* $     is independent of $j$.  
Thus, using a variant of the polar decomposition theorem, we  write  $Y_j = M U_j $  for some 
positive definite operator $ M $ and $U_j $ unitary.    Then 
$\tau(Y_j Y_k^*) =  \tau(M^2 U_j  U_k^* ) = \delta_{jk} $, implies $\tau(M^2) = \tau(I_{\cal N} )$.
Moreover, since $Q_{jk} ^\pm  = M U_j U_k^* M \pm  M U_k U_j^* M $, one can multiply
 by $M^{-1} $ to conclude that  the homogenous equations \eqref{jk=mn} to \eqref{U-cond}  hold 
  for the $U_j $ in $Y_j = M U_j $ even if  $Y_j$ is not unitary.  Moreover, they hold
for  $U_j U_k^* \pm U_k  U_j^* $ if and only if they hold for $Y_j Y_k^* \pm Y_k Y_j^*$, justifying our ambiguous use of $Q_{jk}^\pm $.

Applying the   argument above to  $Y_j^* Y_k $ implies $Y_j^* Y_j $ is independent of $j$ so that
$   U_j^* M^2 U_j = U_k^* M^2 U_k $. This implies that $U_j U_k^* $ commutes with $M^2$ 
for all pairs  $j \neq k $ so that 
\be    \label{s=tb1} 
        28  I_\mathcal{N} & = &  28 \,  U_j^* M^2 U_j  +  6 \sum_{k \neq j } ( U_j^* M^2 U_k +  U_k^* M^2 U_j )     \\
        \label{s=tb2}          & = &  28\,  M^2 + 6   M^2 \sum_{k \neq j } (     U_k U_j^* +  U_j U_k^*   )  \qquad \qquad       j = 1,2,3,4 .
 \ee 
 When $M^2 =  I_\mathcal{N} $, \eqref {s=ta} implies  $\sum_{k \neq j} Q_{jk}^+ = 0 $ and 
 \eqref{s=tb1} implies   $\sum_{k \neq j} R_{jk}^+ = 0 $ for each fixed $j = 1,2,3,4$.
 Combining this with \eqref{jk=mn}  gives \eqref{U+cond}.   When $M^2 \neq   I_\mathcal{N} $, we obtain \eqref{s=tb2}.  \qed  \medskip
 
 Even though it seems unlikely, we can not exclude the possibility that $\Phi$ is factorizable
when  $M^2 \neq  I_\mathcal{N} $.    % if $Y_j $ is not unitary,  then   $ \tau ( U_j  U_k^* ) = \delta_{jk} $  need not hold.
Nevertheless, we can state necessary 
and sufficient conditions for factorizability in terms of unitary operators.
 \begin{remark}   \label{M2}
 The channel $\Phi$ in Remark~\ref{rmk:d4fact} has an exact factorization through $M_d(\mathbf{C}) \otimes \mathcal{N}$ with 
$\bU = \sum_{k= 1}^4 A_k \otimes Y_k $  if and only if there are unitary operators $U_k \in {\cal N} $
such that the following conditions hold:

 {\rm a)}  $\b I_\mathcal{N} +\tfrac{3}{14} \sum_{k \neq j } Q_{jk}^+   $ is positive definite,
   $\tau \big( ( I_\mathcal{N} +\tfrac{3}{14} \sum_{k \neq j } Q_{jk}^+  )^{-1} U_m U_n^* \big) = \delta_{mn} $, and
 
  {\rm b)} both \eqref{jk=mn} and \eqref{U-cond}   hold.
   
  When these conditions hold, one can choose $M^2 = \big( I_\mathcal{N} +\tfrac{3}{14} \sum_{k \neq j } Q_{jk}^+  \big)^{-1} $ 
  and $Y_n= M U_n$.   
  Moreover, \eqref{jk=mn} to \eqref{U-cond} hold if and only if they also hold when $U_j $ is replaced by $Y_j $
  in $R_{jk}^\pm$.
  \end{remark}  

\noindent{\bf Proof:}  Both \eqref{s=ta} and \eqref {s=tb2}  are equivalent to
  $   M^{-2} =  I_\mathcal{N} + \tfrac{ 3}{14} \sum_{k \neq j } (   U_k U_j^* +  U_j U_k^* )$.  
Then $  \tau ( Y_m  Y_n^* ) =  \tau(M^2 U_m U_n^* )  = \delta_{mn} $ gives  (a).     Although \eqref{U+cond} does not hold,
\eqref{jk=mn} implies that $ \sum_{k \neq  j} Q_{jk}^+ $ is independent of $j $ so that the definition of $M^2$ is also.

However,  relating conditions on $Y_j^* Y_k $ to those on $U_j^* U_k $, 
  requires more work.
Using the conditions for $s \neq t$ on $Y_j^* Y_k $, we can conclude that 
$ Y_j^* Y_k + Y_k^* Y_j =  Y_m^* Y_n + Y_n^* Y_m  $
 when $j,k,m,n $  are distinct,  which is equivalent to
\be  \label{jk=mn*}
 U_j^* M^2 U_k + U_k^* M^2  U_j = U_m^* M^2 U_n + U_n^* M^ 2 U_m .
\ee
Now observe that, since $M^2  $  commutes with $  U_k U_\ell^* $,
\be   \label{flip*}
 U_j^* M^2 U_k = U_j^* M^2 U_k U_\ell^* U_ \ell = U_j^*  U_k U_\ell^* M^2  U_ \ell  .
 \ee
Applying \eqref{flip*} to all terms in \eqref{jk=mn*} and then multiplying on the right by $U_\ell^* M^{-2} U_\ell^* $
gives   $ U_j^*   U_k + U_k^*    U_j = U_m^*   U_n + U_n^*  U_m $.  
One can similarly show that \eqref{U-cond} holds for $R_{jk}^- $. \qed

We have found  $U_j  \in M_4({\bf C}) $ which satisfy the conditions for $Q_{jk}^+$.  They
have the form $U_j = U_j^* = 2 E_j - \rmid_4 $ with each $E_j$ a rank two projection from one of the 4 mutually 
unbiased bases \cite{K} for ${\bf C}_4 $ (excluding the standard basis).  However, they do not satisfy \eqref{U-cond}.
 It is plausible that one could use the fact that \eqref{U-cond}  is associated with an irreducible representation of 
 the symmetric group,  as in \eqref{skewx=2}, to find $U_j$ which satisfy \eqref{U-cond}.
  If so, the question of factorizability depends on
whether or not the conditions for $Q_{jk}^+$ and $Q_{jk}^-$ are compatible.  One can use the asymmetric
condition \eqref{asym:jkl} to test this without finding $Q_{jk}^+$ or $Q_{jk}^-$.

 \begin{remark}
 The conditions  in Remark~\ref{rmk:d4fact} can not hold for $U_k \in M_\nu({\bf C} ) $ unless $\nu = 3 n$ is a multiple of $3$.
 This implies that the channel  $\Phi$  in Remark~\ref{rmk:d4fact} can not have an exact factorization through 
 $M_d(\mathbf{C}) \otimes M_\nu({\bf C} )$ when $\nu \neq 3n$.
  \end{remark}
\noindent{\bf Proof:}  For $j,k,m = 1,2,3$ multiplying  the first equation in \eqref{asym:jkl} by $U_2$ on the right and the second by $U_3$ on the left gives
\be
    U_1 + U_2 U_3^* U_2 + U_3 U_1^* U_2 = 0   \qquad  and  \qquad   U_3 U_1^* U_2 + U_3 U_2^* U_3 + U_1  = 0
\ee
which implies $U_2 U_3^* U_2  = U_3 U_2^* U_3  $ or, equivalently $(U_3^* U_2 )^2 = U_2^* U_3$.    Multiplying this by $U_3^* U_2 $ gives
$(U_3^* U_2)^3 = \rmid_\nu $.   Thus, we can conclude
\be
   \tr U_3^* U_2 = 0  \qquad   \quad  \tr (U_3^* U_2 )^2 =  \tr U_2^* U_3 = 0   \qquad  \quad  \tr (U_3^* U_2)^3 = \tr \rmid_\nu = \nu
\ee
Let $\lambda_k $ be the eigenvalues of $U_3^*  U_2 $ which is unitary.  Then we have
\be
    |\lambda_k | = 1 ~~ \forall ~ k   \qquad  \quad \sum_{k=1}^\nu\lambda_k = 0   \qquad \quad \sum_{k=1}^\nu\lambda_k^2 = 0   \qquad \quad  \sum_{k=1}^\nu\lambda_k^3 = \nu.
\ee 
The last equation implies $\lambda_k^3 = 1 ~~ \forall ~  k $.  Thus, the only possible choices for  $\lambda_k$ are $1, e^{\pm 2 \pi i /3} $.
Then  $ \sum_k \lambda_k =  \sum_k \lambda_k^2 = 0 $  holds if and only if   the eigenvalues of  $ U_3^*U_2$ are
$ 1, e^{2 \pi i /3}, e^{-2 \pi i /3}   $,  each  with the same multiplicity. Thus, $\nu$ must be a multiple of $3$.
  \qed  % This also holds for any pair $U_j^* U_k $. \qed

This does not resolve the question of factorizability in the case $d = 4, t = - \third  $, although it seems unlikely.
We have not excluded the possibility of satisfying \eqref{jk=mn} to \eqref{U-cond} with unitary matrices  $U_k \in M_{3n}({\bf C}) $.
Moreover, we have not excluded the possibility of satisfying the conditions in Remark \ref{M2} when $M^2 \neq  I_{\cal. N}$.   Nor have we excluded the 
possibility  that $\Phi$ does  have an exact factorization through 
 $M_d(\mathbf{C}) \otimes \mathcal{N}$ with ${\cal N} $ not a matrix algebra, as in \cite{MR}.

  \pagebreak

 \begin{remark}
When $ d = 4 $, the channel $\Phi$ is  factorizable at most when $t = \pm 1$ or $t = -\third$.  Nevertheless, the channel $\Phi \circ \Phi$ is
 factorizable for all $t \in [-1,1]$.  
   \end{remark}
 Since $\Phi = \Phi^*$, this follows immediately from \cite[Remark~5.6]{HM2011} as 
 observed in Appendix~\ref{app:choi4}.
  
\bigskip
  
  \noindent{\bf   Acknowledgment:}     It is a pleasure to thank Matthias Christandl,  Mikael R{\o}rdam  and Anders Thorup for helpful discussions, 
  and the graduate student Jon Lindegaard Holmberg for skillfully performing insightful numerical work. 
  The second named author was supported by a grant from The Independent Research Fund Denmark (FNU).

% \bigskip  %  \pagebreak

 \appendix
 
 \section{Factorizability  }  \label{app:fact}
 
\subsection{Conditions for factorizability}    \label{app:factcond}

Let $\Phi \colon M_d(\mathbf{C}) \to M_d(\mathbf{C})$ be a UCPT map  and $\mathcal{N}$ a von Neumann algebra with a faithful trace $\tau$
normalized so that   $\tau(I_{\mathcal{N}}) = 1$.   Following [6],  we say  (as noted in the introduction) 
that $\Phi$ has an exact factorization through $M_d(\mathbf{C}) \otimes {\mathcal{N}}$ 
if there is a unitary $\mathbf{U} \in M_d(\mathbf{C}) \otimes \mathcal{N}$ such that for all $\rho \in M_d(\mathbf{C})$  \eqref{exactfact} holds, i.e.,
 $\Phi(\rho) = (\mathcal{I} \otimes \tau)\mathbf{U}^*(\rho \otimes I_{\mathcal{N}}) {\mathbf{U}}$.
Furthermore, $\Phi$ is called \emph{factorizable} if it has an exact factorization through 
$M_d(\mathbf{C}) \otimes \mathcal{N}$ for \emph{some} $(\mathcal{N}, \tau)$.

When $\mathcal{N} = M_\nu(\mathbf{C})$ is a matrix algebra, \eqref{exactfact} becomes
  $\Phi(\rho) = (\mathcal{I} \otimes \mathrm{Tr})\mathbf{U}^*(\rho \otimes \frac{1}{\nu}  \rmid_\nu) {\mathbf{U}}$.
  Following standard practice used elsewhere in this paper, when   ${\cal N} = M_d(\mathbf{C}) $, we 
  denote the identity by $\rmid_d$ rather than the awkward 
  $\rmid_{M_d(\mathbf{C})} $.  We also follow the standard convention for type I algebras that $\tr$
  is normalized  so that $\tr \proj{v} = 1 $ when $v$ with $\norm{v}^2 = \bra v, v \ket = 1$.

The following useful result follows   from the   factorizability criteria in  [5, Theorem 2.2].
\begin{prop}     \label{prop:factcond}
Necessary and sufficient conditions for a UCPT map $\Phi$ on $M_d(\mathbf{C})$ of the form $\Phi(\rho) = \sum_{k=1}^\kappa A_k^* \rho A_k$, with $\{A_1, A_2, \dots, A_\kappa\}$ linearly independent to have an exact factorization through $M_d(\mathbf{C}) \otimes \mathcal{N}$ are that there are $Y_1,Y_2, \dots, Y_\kappa \in \mathcal{N}$ such that 
$\tau( Y_j Y_k^*) = \delta_{jk} $ and (with $\{|e_s\rangle\}_s$ the standard basis  for $\mathbf{C}_d$) the following conditions hold for all $1 \le s,t \le d$:
\begin{equation}   \label{factcond} 
\sum_{j,k} \langle e_s, A_jA_k^* e_t \rangle Y_jY_k^* = \delta_{st} \, I_\mathcal{N}, \qquad
\sum_{j,k} \langle e_s, A_j^*A_k e_t \rangle Y_j^*Y_k = \delta_{st} \, I_\mathcal{N}.
\end{equation}
\end{prop}

\noindent{\bf Proof:} 
By [5, Theorem 2.2], if $\Phi$ has an exact factorization through $M_d(\mathbf{C}) \otimes {\mathcal{N}}$, then one can find $Y_1,Y_2, \dots, Y_\kappa$ in $\mathcal{N}$ such that $\mathbf{U} := \sum_{k=1}^\kappa A_k \otimes Y_k \in M_d(\mathbf{C}) \otimes \mathcal{N}$ is unitary. Hence
\begin{equation}    \label{unitcond}
\mathbf{U}\mathbf{U}^* = I_d \otimes I_{\mathcal{N}} = \sum_{j,k} A_jA_k^* \otimes Y_jY_k^*, \qquad
\mathbf{U}^*\mathbf{U} = I_d \otimes I_{\mathcal{N}} = \sum_{j,k} A_j^*A_k \otimes Y_j^*Y_k.
\end{equation}
Now observe that for any $A \in M_d(\mathbf{C})$ and $Y \in \mathcal{N}$ one has
$$A \otimes Y = \sum_{s,t} \langle e_s, A e_t \rangle \, | e_s \rangle \langle e_t |  \otimes Y = \sum_{s,t} | e_s \rangle \langle e_t | \otimes \langle e_s, A e_t \rangle \,  Y.$$
Applying this to both sides of the equations in \eqref{unitcond} and using the fact that $\{ |e_s \rangle \langle e_t |\}_{s,t}$ is a basis for $M_d(\mathbf{C})$, gives the equations in \eqref{factcond}. 
\qed

\subsection{Pairs of factorizable maps}   \label{app:dual} 
 \begin{defn}  \label{def:dual}
Let  $ d =  p\cdot q $ be a product of primes and $\bU \in  M_d({\bf C}) \simeq  M_p({\bf C}) \ot M_q({\bf C})$   unitary. 
Then $\Phi(\rho) = ( \id \ot {\rm Tr} )\bU^* \big(\rho \ot \tfrac{1}{q}  \rmid_q \big) \bU$ and 
 $\Psi(\gamma) = ( {\rm Tr} \ot \id )\bU^* \big( \tfrac{1}{p } \rmid_p \ot \gamma \big) \bU$
are UCPT maps on $M_p ({\bf C})$  and  $M_q ({\bf C}) $, respectively. air   The  channels $\Phi, \Psi $ are said to be {\em dual channels
associated with the unitary} $\bU $.
 \end{defn} 
 
 Theorem~\ref{thm:duald3}  essentially says that  
 the channels $\Phi, \Psi$ in Section~\ref{sect:d=3} are dual channels associated with
 the unitary  $\bU$ given by \eqref{Udef} and $p = q = 3$.   
 
 Let $\{ U_j \}$ be a set of $q$ unitary matrices in $ M_p({\bf C})$ and
   ${\bf X}   = \bigoplus_{j=1}^q U_j  = \sum_{j=1}  U_j \ot \proj{e_k} $.  \linebreak Then 
   ${\bf X} \in    M_p({\bf C}) \ot  M_q ({\bf C})  $  is unitary and, as
 essentially observed in \cite[Proposition 2.8]{HM2011}, the associated pair of dual channels
 are $\Phi(\rho) =   \tfrac{1}{q} \sum_{j=1}^q U_j^* \rho \, U_j $ and  $\Upsilon(\gamma) = C \circledast \gamma$, 
 where $\circledast$ denotes the Schur or Hadamard product and $C  \in  M_q ({\bf C}) $ is the matrix with elements 
 $c_{jk} = \tfrac{1}{p} \tr U_j^* U_k $.  When the $U_j$ are also orthogonal 
 so that $\tr U_j^* U_k = p  \, \delta_{jk}  $, then
$\Upsilon(\gamma) =    \rmid_q \circledast \gamma  =  \sum_{j=1}^q \gamma_{jj} \proj{e_j} $ is diagonal.

 Whenever a UCPT map $\Phi$  has an exact factorization through $ M_d({\bf C}) \ot M_\nu({\bf C}) $,
 there is a unitary $\bU \in M_d({\bf C}) \ot M_\nu({\bf C}) $ for which 
 $\Phi(\rho) = (\mathcal{I} \otimes \mathrm{Tr})\mathbf{U}^*(\rho \otimes \frac{1}{\nu}  \rmid_\nu) {\mathbf{U}}$ so that 
 there is another UCPT map $\Psi$  such that $\Phi, \Psi$  are the dual pair
 associated with that unitary.    However, a UCPT map can have an exact factorization   
 in more than one way.  An example is given by the channel $\Psi $ in Section~\ref{sect:d=3} which corresponds to
 $t = 1$ in \eqref{keydef}.    This channel has an exact factorization through   $ M_3({\bf C}) \ot M_3({\bf C}) $
 with the unitary   $\bW $ given by \eqref{Wdef} and another with ${\bf X}   = \sum_{k = 1}^3  B_j \ot \proj{e_j} $ 
 as in Remark~\ref{rmk:factpm1}.  The  dual pair associated with $\bW$ is $\Psi, \Phi $ as in Section~\ref{sect:d=3};
 the dual pair associated with ${\bf X} $ is $\Psi, \Upsilon$ where $\Upsilon(\gamma) =\rmid_3 \circledast \gamma $
as above.

One can extend this to situations in which $d$  is a product of more than two primes.
However, each way of writing $d$ as a product of two integers will give a different pair of dual channels.  
Finally, we note that when $p = q$, it is possible to have  a self-dual channel.   For example, when
 $\bU =\1rt2 \big( \sigma_x \ot \sigma_z + i \sigma_z \ot \sigma_x ) $, one finds
 $\Phi(\rho) = \Psi(\rho) =  \half \big( \sigma_x \rho \sigma_x +  \sigma_z \rho \sigma_z )$.
 with $\sigma_{x,y.z} $ denoting the usual Pauli matrices.

It is worth noting that   $\bU = \sum_k A_k \ot B_k $ with $A_k \in  M_p({\bf C})$ and $B_k \in M_q({\bf C})$ does not imply
that $\Phi(\rho) = \sum_k A_k^* \rho A_k $ and $ \Psi(\gamma) = \sum_k B_k^* \gamma B_k $.  This only holds if
$\tr B_j^* B_k =  q \, \delta_{jk}  $ in the first case and $\tr A_j^*A_k =p  \, \delta_{jk} $ in the latter.  Thus, for example,
when $\bU$ is given by \eqref{unitaryUCPT},  $\Psi$ in the dual pair $\Phi_{\alpha,\beta} , \Psi_{\alpha,\beta} $ is
{\em not} given by $\Psi(\gamma) = \half \sum_{j,k =1}^2 |e_j \kb e_k | \gamma |e_k \kb e_j | = \half  ( \tr \gamma ) \rmid_2 $, but by
\bee
\lefteqn{     \Psi_{\alpha,\beta} =   \tfrac{|\alpha|^2 + 1}{3} \big( \proj{e_1} \gamma \proj{e_1}   + \proj{e_2} \gamma \proj{e_2}  \big) } \quad \\
       & ~ & \quad 
                +  ~  \tfrac{|\beta|^2 + 1}{3}  \big( |e_1 \kb e_2 | \gamma |e_2 \kb e_1 | + |e_2 \kb e_1 | \gamma |e_1 \kb e_2 | \big) \\
                &= & \tfrac{1}{3} (\tr \gamma) \rmid_2
                         +\tfrac{1}{3} \pmx a^2 \gamma_{11} + b^2 \gamma_{22} & 0 \\ 0 & b^2 \gamma_{11} + a^2 \gamma_{22} \emx  .
 \eee

It should be emphasized that this duality is not equivalent to the notion of a pair of  ``complementary channels''
 used in the quantum information literature \cite{DS,H,KMNR} , which is defined in terms of the Stinespring
 representation, and  goes back to Arveson   \cite[Section 1.3]{A}  who used the term ``lifting''.  In that
 case, the auxiliary space is interpreted as the environment and the complementary channel maps the
 input state $\rho$ to a state for the environment.  For a channel   $\Phi : {\bf C}_{d_A} \mapsto  {\bf C}_{d_B} $
  of the form  $\Phi(\rho_A) = \sum_{k = 1}^{d_E}  A_k^* \rho A_k $, one can regard the Stinespring representation
  as mapping  $\rho_A \mapsto \rho_{BE} = \sum_{jk} A_j ^*\rho_A A_k \ot  |e_j \kb e_k |$ with $|e_j \ket $
  the standard basis for ${\bf C}_{d_E}$.  Then 
  $\Phi(\rho_A) = \trp_E \, \rho_{BE} $ and the complementary channel   $\Phi^C :  {\bf C}_{d_A} \mapsto  {\bf C}_{d_E}$
  is  
 \bee
     \Phi^C(\rho_A) =  \trp_B \, \rho_{BE}  = \sum_{jk}   \tr(A_j ^*\rho_A A_k )   |e_j \kb e_k | .
 \eee
Even when $d_A = d_B$ these concepts are quite different.   Complementary channels are defined
 with the implicit assumption of a pure ancilla rather than a maxiamlly mixed ancilla.  In the notion of dual pairs 
introduced above, the roles of the input and environment are interchanged in terms of both the
subspace over which the trace is taken {\em and}  the space in which the (maximally mixed) ancilla resides.

\subsection{Factorizability of $\Phi \circ \Phi^*$  with Choi rank $ \leq 4$.}   \label{app:choi4}

When a UCPT map  $\Phi $ is factorizable,  the adjoint $\Phi^*$ is also factorizable.
Moreover, the maps  $\Phi \circ \Phi$,  $\Phi^* \circ \Phi$, and $\Phi \circ \Phi^*$ are also factorizable.
However, there are some special circumstances in which $\Phi$ is {\em not } factorizable, but $\Phi^* \circ \Phi$
 and $\Phi \circ \Phi^* $ are factorizable.  This includes the Arveson-Ohno channel  \eqref{ohnomap} and 
 the channels in Section~\ref{subsect:key} for $d = 4$  and all
 $ t \in (-1, 1) $.

The next result is a straightforward generalization of  \cite[Remark~5.6]{HM2011}.
  It follows from \cite[Lemma 5.5]{HM2011} that the condition Choi rank $\leq 4$ is critical.
\begin{prop}   \label{prop:choi4}
Let $\Phi:  M_d({\bf C})  \mapsto  M_d({\bf C})  $ be a UCPT map with Choi rank $\leq 4$.  
Then the maps $\Phi \circ \Phi^* $  and $\Phi^* \circ \Phi $ each have exact factorizations through
 $ M_d({\bf C})  \ot  M_4({\bf C})  $.
\end{prop}
\noindent{\bf Proof:}  Let $\{A_k \in M_d({\bf C} ) : k =1, 2, 3, 4 \} $ satisfy $\sum_{k = 1}^4 A_k^*A_k = \sum_{k = 1}^4 A_k A_k^* = \rmid_d$
and define the UCPT map $\Phi(\rho) =  \sum_{k = 1}^4 A_k^* \rho A_k $.  As observed in \cite{HM2011}, one can always choose some $A_k  = 0 $
so that the result follows if  there is a unitary map $\bU \in M_{4d}({\bf C}) $ such
that  
\be  \label{d5fact}
(\id \ot \tr)  \bU^* ( \rho \ot \tfrac{1}{4} \rmid_4 ) \bU = (\Phi \circ \Phi)^*(\rho) =  \sum_{j = 1}^4  \sum_{k = 1}^4 A_j^* A_k \rho A_k^* A_j .
\ee
Following the strategy in  \cite[Remark~5.6]{HM2011}, define $   \bU = \sum_{j,k = 1}^4  A_j^* A_k \ot (2 | e_j \kb e_k |  - \delta_{jk} \rmid_4 )$. 
Then, by repeatedly using    $ \sum_{k = 1}^4 A_k^*A_k = \sum_{k = 1}^4 A_k A_k^* = \rmid_d$ one finds  %   \pagebreak
\bee
    \bU \bU^* & = &   \sum_{j,k = 1}^4   \sum_{m,n = 1}^4 A_j^* A_k A_n^* A_m \ot 
         ( 2 | e_j \kb e_k |  -\delta_{jk} \rmid_4 ) (2  | e_n \kb e_m |  - \delta_{mn} \rmid_4 ) \\
        & = & 4 \sum_{jkm} A_j^* A_k A_k^* A_m \ot  |e_j \kb e_m |  - 2 \sum_{j,m,n} A_j^* A_j A_n^* A_m \ot   | e_n \kb e_m |   \\
          & ~ & \quad -  2\sum_{jkm}  A_j^* A_k A_m^* A_m \ot  | e_j \kb e_k |  +  \sum_{k,n} A_k^* A_k A_n^* A_n \ot \rmid_4   \\   
          & = & 4 \sum_{jm} A_j^* A_m \ot  |e_j \kb e_m | - 2 \sum_{m,n}   A_n^* A_m \ot   | e_n \kb e_m | 
                        - 2 \sum_{jk}  A_j^* A_k \ot  | e_j \kb e_k | + \rmid_d \ot \rmid_4 \\
          & = &         \rmid_d \ot \rmid_4       
\eee 
so that $\bU$ is unitary.  Similary,  one finds  (using $ \tr |e_j \kb e_k | = \delta_{jk} $ and $\tr \rmid_4 = 4 $ )
\bee
     (\id \ot \tr)  \bU^* ( \rho \ot \tfrac{1}{4} \rmid_4 ) \bU %  = } \\
    & = &   \tfrac{1}{4}  \sum_{jkmn} A_j^* A_k \rho A_n^* A_m \big( 4 \delta_{jm} \delta_{kn} -2 \delta_{jk} \delta_{mn} 
                 -2  \delta_{jk} \delta_{mn}    +  4  \delta_{jk} \delta_{mn}  \big)  \\
             & = &  \sum_{jk} A_j^* A_k \rho A_k^* A_j   =  (\Phi \circ \Phi)^*(\rho) \hskip4.5cm  \hbox{QED}              \eee  

\section{Linear dependence of $\{ A_m^* A_n \} $ vs $\{ A_m A_n^* \} $.}  \label{app:B}

We give an explicit example to show that  $\{ A_m^* A_n \} $ can be linearly independent, but $\{ A_m A_n^* \} $
  linearly dependent. 
Let  $d = 4$ and $W = \tfrac{1}{21} \pmx   8 & -11 &  16 \\  -19 & -8 &   4 \\  -4 & 16 & 13 \emx $.
When $A_m $ is constructed as in \eqref{tform} with all $V_m = W$, we found that
\begin{itemize}

\item   $\{ A_m^* A_n  \}  $  is linearly independent unless $t = \pm 1,~ t = \tfrac{-13}{3} , ~ t = \tfrac{-59}{84} ,~  t = \tfrac{19}{21}$ , or $ t = \tfrac{107}{21}$.

\item   $\{ A_m A_n^*   \} $  is linearly independent unless $t = \pm 1 ,~ t = \tfrac{-59}{84} ,~  t = \tfrac{-1}{7} , t=  \tfrac{19}{21}$ , or $ t = \tfrac{107}{21}$.

\end{itemize} 
Thus, when $t = -\tfrac{1}{7}$,   $\{ A_m^*  A_n^*   \} $  is linearly independent  but $\{ A_m A_n^*  \}  $  is linearly dependent.
When $t =  - \tfrac{13}{3} $,  $\{ A_m A_n^*   \} $  is linearly independent  but $\{ A_m^* A_n  \}  $  is linearly dependent.
At $t = \pm 1,  ~t = \tfrac{-59}{84} ,~  t = \tfrac{19}{21}, ~ t = \tfrac{107}{21}$ both sets are linearly dependent.  For all other values of $t \in {\bf R} $
both sets are linearly independent.

This result might seem counter-intuitive because the cyclicity of the trace implies 
 $$  \tr A_j^* A_k (A_m^* A_n )^* = \tr A_m  A_j^* A_k A_n^*  = \tr A_m  A_j^* (A_n A_k^*)^* $$
so that the Gram matrices for the two sets have the same elements, albeit arranged differently.
Let $G$ and $H$ denote these Gram matrices  % for  $\{ A_m^* A_n   \} $ and $ \{ A_m A_n^*   \} $
%respectively, 
and consider the elements 
  $ g_{11,kk} = \tr A_1^* A_1 (A_k^* A_k )^*    = \tr  A_k A_1^* (A_k A_1^* )^* = h_{k1,k1}$.
  Since   $ g_{11,kk}$  all lie in the first row of $G$, $\det G$ will not contain
  any terms with  $ g_{11,kk}  \cdot g_{11,jj} $ when $ j \neq k $.  However, since $ h_{k1,k1}$ lies on
  the diagonal of $H$,  
  $ \det H$ will contain a term which includes $\prod_{k = 2}^d h_{k1,k1}   = \prod_{k = 2}^d g_{11,kk}$.
  
We conjecture that  if  $\{ A_m A_n^*  \}  $  is linearly dependent for all $t$,  then
       $\{ A_m A_n^*   \} $  should also be linearly dependent for all $t$.


\bigskip      

 \begin{thebibliography}{~~}
 
 \bibitem{A-D}  C.~Anantharaman-Delaroche,  ``On ergodic theorems for free group actions on
  noncommutative spaces'', {\em Probab. Theory Related Fields} \textbf{135},  520--546 (2006).
 
 \bibitem{A} W. Arveson, ``Subalgebras of $C^*$-algebras''
 {\em Acta Mathematica} {\bf 123}, 141-224 (1969).
 
  \bibitem{eof}    C.H. Bennett, D.P. DiVincenzo, J. Smolin, and W.K. Wootters,
  ``Mixed State Entanglement and Quantum Error Correction''
{\em Phys. Rev. A }  {\bf 54}, 3824 (1996.)
 
 \bibitem{Choi}  M-D. Choi,
``Completely positive linear maps on complex matrices''
    {\em  Lin. Alg. Appl.}    {\bf 10}, 285--290
(1975).

\bibitem{DS} I. Devetak and P. Shor,
 ``The capacity of a quantum channel for simultaneous transmission of classical and 
 quantum information''
 {\em Commun. Math. Phys.} {\bf 256}, 287-303 (2005).
            
            \bibitem{HM2011} U. Haagerup and M. Musat,
            ``Factorization and Dilation Problems for Completely
Positive Maps on von Neumann Algebras''
             {\em Commun. Math. Phys.}  {\bf 303}, 555--594 (2011).
                     
             \bibitem{HM2015} U. Haagerup and M. Musat,
             ``An Asymptotic Property of Factorizable Completely
Positive Maps and the Connes Embedding Problem''
{\em Commun. Math. Phys. } {\bf 338},  721- -752 (2015).

\bibitem{H}  A.S. Holevo,  
``On complementary channels and the additivity problem''
{\em Prob. Theory Appl.} {\bf 51},  133-143 (2006).
            
            \bibitem{HJ}  R.A. Horn and C.R. Johnson,
            {\em Topics in Matrix Analysis} (Cambridge Press, 1991).
            
             \bibitem{HHO} M. Horodecki, P. Horodecki, and J. Oppenheim 
             ``Reversible transformations from pure to mixed states and the unique measure of information'' 
             {\em Phys. Rev. A}
{\bf 67}, 062104 (2003).

            
            \bibitem{J} A. Jamiolkowski
           ``Linear transformations which preserve trace and positive semi-definiteness of operators''
  {\em Reports Math. Phys.} {\bf 3}, 275--278 (1972).
  
  \bibitem{JMVWY} 
Z. Ji, A. Natarajan, T. Vidick, J. Wright, and H. Yuen, 
``MIP*=RE`` arXiv:2001.04383v2.

\bibitem{K}   W.K. Kantor,  ``MUBs inequivalence and affine planes''
{\em J. Math. Physics}  {\bf 53}, 032204 (2012).


\bibitem{KMNR} C. King, K. Matsumoto, M. Nathanson, and M.B. Ruskai,
``Properties of conjugate  channels with applications to additvity and multiplicity''
{\em Markov Processes Related Fields} {\bf 13}, 391--423 (2007).

            
                     \bibitem{K} B. K\"{u}mmerer 
                     ``Markov dilations on the {$2\times 2$} matrices''  
                    {\em  Operator Algebras and their Connections with Topology and Ergodic Theory} 
             (Proceedings, 1983)        
              Ed. by H. Araki, C.C. Moore,  S. Stratila, and D. Voiculescu.  
                 {\em   Lecture Notes in Math. } {\bf{1132}},  312--323  (Springer, 1985).
                     

            
            \bibitem{LS} L.J. Landau and R.F. Streater,  
            ``On Birkhoff's theorem for doubly stochastic
completely positive maps of matrix algebras''
    {\em  Lin. Alg. Appl.} {\bf 193}, 107--127 (1993).
    
            
            \bibitem{MW}   C.B. Mendl and M.M. Wolf
``Unital Quantum Channels - Convex Structure and Revivals of Birkhoff's Theorem''
{\em Commun. Math. Phys.}  {\bf  289}, pp. 1057--1096 (2009).

      
           \bibitem{MKZ}{ M. Musz, M. Kus, and K. Zyczkowski}
       ``Unitary quantum gates, perfect entanglers, and unistochastic maps''
       {\em Phys. Rev. A}  {\bf 87}, 022111 (2013).


\bibitem{MR}  M. Musat and  M. Rørdam
``Non-closure of quantum correlation matrices and factorizable channels that require infinite dimensional ancilla''
{\em Commun. Math. Phys.}  {\bf 375}, 1761--1176 (2020). 
      
            \bibitem{Ohno} H. Ohno
            ``Maximal rank of extremal marginal tracial states''
          {\em J. Math. Phys.}  {\bf 51}, 092101 (2010).   
          
            \bibitem{Rusk} M.B. Ruskai, 
            ``Some Open Problems in Quantum Information Theory''
      arXiv:0708.1902.

           \bibitem{SM}  J. Scharlau and M.P. M\"uller {\em Quantum
Horn lemma, finite heat baths,
and the third law of thermodynamics},{\em  Quantum}  {\bf 2},  54 (2018).
       
       \bibitem{S} S. Sternberg {\em Group Theory and Physics}
       (Cambridge Press, 1994).
       
       \bibitem{WH}  R.F. Werner and  A.S. Holevo
      `` Counterexample to an additivity conjecture for output purity of quantum channels''
      {\em J. Math. Phys.} {\bf 43}, 4353 (2002).
      

           \end{thebibliography}
    \end{document}